\documentclass[12pt, reqno]{amsart}
\usepackage[all]{xy}
\usepackage{amssymb}

\allowdisplaybreaks[2]
\newcommand{\xXi}{\Xi\!\Xi}

\newcommand{\indlp}{\Ind_{L}^{P}}
\newcommand{\indkp}{\Ind_{K}^{P}}
\newcommand{\indgq}{\Ind_{G}^{Q}}
\newcommand{\indkq}{\Ind_{K}^{Q}}
\newcommand{\indlq}{\Ind_{L}^{Q}}
\newcommand{\indkpsq}{\Ind_{K}^{\psq}}
\newcommand{\indgpsq}{\Ind_{G}^{\psq}}
\newcommand{\cc}{C_{c}}
\let\phi\varphi

\def\ZZ{{W}}

\newcommand{\kup}{K\backslash P}
\newcommand{\qmm}{Q/G}
\newcommand{\psq}{P\mathrel{\#}Q}
\newcommand{\pspq}{P\mathrel{\#}_\phi Q}
\newcommand{\sigmat}{\widetilde\sigma}
\newcommand{\taut}{\widetilde\tau}
\newcommand{\atensor}{\odot}
\newtheorem{thm}{Theorem}[section]
\newtheorem{cor}[thm]{Corollary}
\newtheorem{lem}[thm]{Lemma}
\newtheorem{prop}[thm]{Proposition}
\theoremstyle{definition}

\theoremstyle{remark}
\newtheorem{rem}[thm]{Remark}

\def\supp{\operatorname{supp}}
\def\Ind{\operatorname{Ind}}
\def\Aut{\operatorname{Aut}}

\def\Res{\operatorname{Res}}

\newcommand{\under}{\backslash}

\def\clsp{\overline{\operatorname{span}}}

\def\id{\operatorname{id}}
\def\rt{{\operatorname{rt}}} % "right translation"
\def\lt{{\operatorname{lt}}} % "left translation"
\def\rrt{{\operatorname{r}}} % "right"
\def\llt{{\operatorname{l}}} % "left"

\def\H{\mathcal{H}}
\def\L{\mathcal{L}}
\def\K{\mathcal{K}}

\begin{document}

\numberwithin{equation}{section}

\title[Induction in stages for maximal coactions]
{\boldmath{Induction in stages for crossed products 
of $C^*$-algebras by maximal coactions}}

\author[an Huef]{Astrid an Huef}
\address{School of Mathematics and Statistics\\
The University of New South Wales\\
NSW 2052\\
Australia}
\email{astrid@unsw.edu.au}

\author[Kaliszewski]{S. Kaliszewski}
\address{Department of Mathematics and Statistics\\
Arizona State University\\Tempe\\ AZ
85287-1804\\USA} \email{kaliszewski@asu.edu}

\author[Raeburn]{Iain Raeburn}
\address{School of Mathematics and Applied Statistics\\
University of Wollongong\\
NSW 2522\\ Australia}
\email{raeburn@uow.edu.au}

\author[Williams]{Dana P. Williams}
\address{Department of Mathematics\\Dartmouth College\\ Hanover, NH 03755-3551\\USA}
\email{dana.williams@dartmouth.edu}

\thanks{
This research was supported by grants from the Australian Research Council,
the National Science Foundation,
the University of New South Wales
and the Ed Shapiro Fund at Dartmouth College.
}

\subjclass[2000]{46L55}

\date{March 26, 2007}

\begin{abstract}
Let $\delta$ be a maximal coaction of a locally compact group $G$
on a $C^*$-algebra~$B$, and let $N$ and $H$ be closed normal subgroups 
of~$G$
with $N\subseteq H$. We show that the process $\Ind_{G/H}^G$ which uses
Mansfield's bimodule to induce representations of $B\rtimes_\delta G$
from those of $B\rtimes_{\delta|}(G/H)$ is equivalent to the two-stage
induction process $\Ind_{G/N}^G\circ\Ind_{G/H}^{G/N}$.
The proof involves a calculus of symmetric imprimitivity bimodules which
relates the bimodule tensor product to the fibred product of the
underlying spaces.
\end{abstract}

\maketitle

%======================================================================

\section{Introduction}

Induction is a method of constructing representations which is important
in many different situations. The modern $C^*$-algebraic theory of
induction has its roots in Mackey's work on the induced representations
of locally compact groups, which culminated in the Mackey machine for
computing the irreducible unitary representations of a locally compact
group \cite{mackey}, and in Rieffel's recasting of the Mackey machine
in the language of Morita equivalence --- indeed, Rieffel developed his
concept of Morita equivalence for $C^*$-algebras specifically for this
purpose \cite{rie1,rie2}. Takesaki adapted Mackey's construction to the
context of dynamical systems $(A,G,\alpha)$ in which a locally compact
group $G$ acts by automorphisms of a $C^*$-algebra~$A$ \cite{tak}, and
the full strength of the modern theory was achieved when Green applied
Rieffel's ideas to dynamical systems \cite{green}. Takesaki and Green
showed in particular how to induce a covariant representation $(\pi, U)$
of the system $(A,H,\alpha|)$ associated to a closed subgroup $H$ of $G$
to a covariant representation $\Ind_H^G(\pi, U)$ of $(A,G,\alpha)$.

These various theories of induced representations share the following
fundamental properties:

\begin{enumerate}
\item[]\emph{Imprimitivity}:
There is an imprimitivity theorem which characterises the 
representations which are unitarily equivalent to induced representations.

\item[]\emph{Regularity}: 
The representations induced from the trivial subgroup $\{e\}$ 
are precisely the regular representations, up to unitary equivalence.

\item[]\emph{Induction in Stages}: 
If $K$ and $H$ are closed subgroups of $G$ with $K\subseteq H$, 
then $\Ind_H^G\circ \Ind_K^H=\Ind_K^G$, up to unitary equivalence.
\end{enumerate}

Green's formulation of induced representations uses the bijection
$(\tau,V)\mapsto \tau\rtimes V$ between covariant representations
of $(A,G,\alpha)$ and representations of the crossed product
$C^*$-algebra $A\rtimes_\alpha G$, and his induction process is
implemented by (what we now call) a right-Hilbert $(A\rtimes_\alpha
G)$--$(A\rtimes_{\alpha|} H)$ bimodule $X_H^G(\alpha)$: if $(\pi, U)$  is
a covariant representation of  $(A, H,\alpha|)$ on a Hilbert space~$\H$, 
then the induced
representation  $\Ind^H_G(\pi\times U)$ of $A\rtimes_\alpha G$ acts in
$X_G^H(\alpha)\otimes_{A\rtimes_{\alpha|} H}\H$ through the left action of
$A\rtimes_\alpha G$ on $X_H^G(\alpha)$. Green proved that one can fatten
up the left action of $A\rtimes_{\alpha} G$ to an action of $(A\otimes
C_0(G/H))\rtimes_{\alpha\otimes\lt} G$; with this new left action, the
bimodule becomes a Morita equivalence. The resulting imprimitivity theorem
says that a representation $(\tau, V)$ of $(A,G,\alpha)$ on some Hilbert
space $\H_0$ is induced from a representation of $(A,H,\alpha|)$
if and only if there is a representation $\mu$ of
$C_0(G/H)$ on $\H_0$ which commutes with $\tau(A)$ and gives a covariant
representation $(\mu,V)$ for the action $\lt$ of $G$ by left translation
on $C_0(G/H)$ \cite[Theorem~6]{green}. The general theory of Hilbert
bimodules guarantees that the induction process has good functorial
properties, and Green proved induction-in-stages by constructing a
bimodule isomorphism of $X_H^G(\alpha)\otimes_{A\rtimes_{\alpha|}
H}X_K^H(\alpha|)$ onto $X_K^G(\alpha)$ \cite[Proposition~8]{green}.

In nonabelian duality, one works  with coactions of locally compact
groups on $C^*$-algebras: the motivating example is the dual coaction
$\hat\alpha$ of $G$ on a crossed product $A\rtimes_\alpha G$, from
which one can recover a system Morita equivalent to $(A,G,\alpha)$ by taking a
second crossed product $(A\rtimes_\alpha G)\rtimes_{\hat\alpha} G$.
The crossed product $B\rtimes_\delta G$ of a $C^*$-algebra $B$ by a
coaction $\delta$ of $G$ on $B$ is universal for a class of covariant
representations $(\pi,\mu)$ consisting of compatible representations
of $B$ and $C_0(G)$ on the same Hilbert space. Induced representations
of crossed products by coactions were first constructed by Mansfield
\cite{man}, who associated to each closed normal amenable subgroup $N$
a right-Hilbert $(B\rtimes_\delta G)$--$(B\rtimes_{\delta|} (G/N))$
bimodule, and thereby plugged into Rieffel's general framework. Mansfield
checked that inducing from $B\rtimes_{\delta|}(G/G)=B$ gave the generally
accepted class of regular representations \cite[Proposition~21]{man},
and proved an elegant imprimitivity theorem: a representation
$\tau$ of $B\rtimes_\delta G$ is induced from a representation
of $B\rtimes_{\delta|} (G/N)$ if and only if there is a unitary
representation $V$ of $N$ such that $(\tau,V)$ is covariant for the
dual action $\hat\delta|$ of $N$. Induction-in-stages was later proved in
\cite[Corollary~4.2]{KQR-DR}.

The hypothesis of amenability appears in Mansfield's theory because
his construction is intrinsically spatial, and the Morita equivalence
underlying his imprimitivity theorem involves the reduced crossed product
$(B\rtimes_\delta G)\rtimes_{\hat\delta,r}N$. Subsequent authors have
shown how to lift the amenability and normality hypotheses \cite{KQ-IF,
aHR-mansfield}, but the resulting imprimitivity theorems still use
the reduced crossed product by the dual action, and are therefore not
well-suited to applications involving covariant representations. In
an effort to produce a theory which is more friendly to full crossed
products by actions, Echterhoff, Kaliszweski and Quigg have proposed the
study of \emph{maximal} coactions \cite{EKQ-MC}, which include
the dual coactions and certain other coactions constructed from them~\cite[\S
7]{KQ-MI}.

Kaliszewski and Quigg have recently shown that for a maximal
coaction $\delta$ of $G$ on a $C^*$-algebra $B$ and any closed normal
subgroup $N$ of $G$, the crossed product $B\rtimes_{\delta|}(G/N)$
by the restriction of $\delta$ is Morita equivalent, via a 
\emph{Mansfield bimodule}
we will denote by $Y_{G/N}^G(\delta)$, to the full crossed product
$(B\rtimes_\delta G)\rtimes_{\hat\delta|}N$. Dropping the left action of
$N$ on their Morita equivalence gives a right-Hilbert $(B\rtimes_\delta
G)$--$(B\rtimes_{\delta|}(G/N))$ bimodule which can be used to define
induced representations $\Ind_{G/N}^G(\pi\rtimes\mu)$, and Theorem~5.3
of \cite{KQ-MI} gives an imprimitivity theorem for this induction process.
Our goal in this paper is to prove 
regularity and induction-in-stages 
for this induction process of Kaliszewski and Quigg.
Regularity is straightforward, and is addressed in the short
Section~\ref{reg}.  
Proving induction-in-stages
--- the assertion that $\Ind^G_{G/H}$ is
equivalent to $\Ind_{G/N}^G\circ\Ind_{G/H}^{G/N}$ --- 
occupies most of the rest of the paper.
Specifically, we will prove:

\begin{thm}
\label{main-thm}
Let $\delta:B\to M(B\otimes C^*(G))$ be a maximal coaction of a locally compact group $G$ on a $C^*$-algebra $B$. Also let
$N$ and $H$ be closed normal subgroups of $G$ with
$N\subseteq H$.  Then the following diagram of right-Hilbert
bimodules commutes:
\begin{equation}\label{eq-mainthm}
\xymatrix{
B\rtimes_\delta G
\ar[dr]_{Y_{G/N}^G(\delta)}
\ar[rr]^{Y_{G/H}^G(\delta)}
&
&
B\rtimes_{\delta|}(G/H).
\\
&
B\rtimes_{\delta|}(G/N)
\ar[ur]_{Y_{G/H}^{G/N}(\delta|)}
&
\\
}
\end{equation}
Equivalently,
\[
Y_{G/H}^G(\delta)
\cong Y_{G/N}^G(\delta)
\otimes_{B\rtimes_{\delta|}(G/N)}Y_{G/H}^{G/N}(\delta|)
\]
as right-Hilbert $(B\rtimes_\delta G)$--$(B\rtimes_{\delta|}(G/H))$  bimodules.
\end{thm}

Here, both $Y_{G/N}^G(\delta)$ and $Y_{G/H}^G(\delta)$ are
Mansfield bimodules
defined using the  coaction $\delta$ of $G$ on $B$. The bimodule
$Y_{G/H}^{G/N}(\delta|)$  is defined using the restricted coaction
$\delta|_{G/N}$  of $G/N$ on $B$ and the normal
subgroup  $H/N \subseteq G/N$, and we have identified the quotient
$(G/N)/(H/N)$ with $G/H$.

%--

The Mansfield bimodule is defined in \cite{KQ-MI} as a tensor product of
three other bimodules (see Remark~\ref{KQ-defn}); thus 
proving that~\eqref{eq-mainthm} commutes using first principles
would involve gluing a different commutative square 
onto each of the arrows in \eqref{eq-mainthm}, 
and then proving that the resulting outer figure ---
which would involve a terrifying nine bimodules --- commutes.
More importantly, this approach obscures the fundamental idea
behind the definition of the $Y$'s, which is 
to pass to the second-dual coaction 
(maximal coactions are precisely those for which full crossed-product
duality holds)
and then invoke 
the symmetric imprimitivity theorem of \cite{rae-sit}.

%--

Thus,
our general strategy for the proof of Theorem~\ref{main-thm} will appeal
to this underlying idea rather than the definition itself.
We will use the naturality of the
Mansfield bimodules (see \cite{BE} for the technical meaning of this)
to reduce to the case where $\delta$ is a dual coaction.  
If $\delta=\hat\alpha$ is a dual coaction, 
it is known (\cite[Proposition~6.5]{KQ-MI}) 
that the Mansfield bimodules 
$Y_{G/N}^G(\hat\alpha)$ and
$Y_{G/H}^G(\hat\alpha)$
appearing in \eqref{eq-mainthm}
can be replaced by bimodules
$Z_{G/N}^G(\alpha)$ and
$Z_{G/H}^G(\alpha)$
constructed using the symmetric imprimitivity theorem.
In Theorem~\ref{thm-YZ}, we extend this result by showing that
$Y_{G/H}^{G/N}(\hat\alpha|)$
is isomorphic to the symmetric imprimitivity bimodule
constructed in 
\cite[Proposition~3.3]{aHKRW-EP},
which we denote by $Z_{G/H}^{G/N}(\alpha)$.  
Combining various 
results from the literature gives an analog of \eqref{eq-mainthm}
for the $Z$'s; and then we can assemble all of our intermediate results 
in Section~\ref{main-proof} to  
complete the proof of~Theorem~\ref{main-thm}.

Because the restriction $\hat\alpha|_{G/N}$ need not be the 
dual of an action of $G/N$, the isomorphism 
$Y_{G/H}^{G/N}(\hat\alpha|)\cong Z_{G/H}^{G/N}(\alpha)$
is not simply
another application of \cite[Proposition~6.5]{KQ-MI};
indeed, establishing this result occupies most of the present paper.
Rather than dealing directly with the definition of the Mansfield
bimodule, 
we appeal to \cite[Corollary~6.4]{KQ-MI}, which shows that
$Y_{G/H}^{G/N}(\hat\alpha|)$ 
can be ``factored'' into a tensor product involving Green
and Katayama imprimitivity bimodules.
The desired isomorphism follows when we show (Theorem~\ref{thm-KXZ})
that 
$Z_{G/H}^{G/N}(\alpha)$ can be factored the same way.
The preparation for the proof of Theorem~\ref{thm-KXZ}
involves identifying each of the three imprimitivity 
bimodules in question with a bimodule constructed from the symmetric
imprimitivity theorem; this is carried out in Section~\ref{previous}.
The proof itself occupies Section~\ref{next}, where we
apply a calculus, developed in Section~\ref{PhashQ},
which allows the tensor product
of such bimodules to be studied at the level of the spaces 
from which they were constructed.

We expect that this calculus will be of independent
interest in the future.  To further illustrate its utility, in 
Section~\ref{another} we apply it
to  the balanced tensor product of two one-sided versions of the symmetric
imprimitivity, thus recovering the isomorphism of the tensor product and
the symmetric version from \cite[Lemma~4.8]{hrw} on the level of spaces.

\subsection*{Notation and conventions}\label{notn}

Our reference for the theory of crossed products by actions and
coactions is \cite{BE}.
We follow the conventions of \cite{KQ-MI} for coactions; in particular, 
all our coactions are non-degenerate and maximal.

We write $\lambda$ and $\rho$ for the left and right regular
representations, respectively, of a group $G$ on $L^2(G)$. If $N$
is a normal subgroup of $G$ we write  $\lambda^{G/N}$ for the
quasi-regular representation of $G$ on $L^2(G/N)$  and $M$ or $M^{G/N}$
for the representation of $C_0(G/N)$ on $L^2(G/N)$ by multiplication
operators, so that $(\lambda^{G/N}_r\xi)(sN)=\xi(r^{-1}sN)$ and
$M(f)\xi(sN)=f(sN)\xi(sN)$ for  $\xi\in L^2(G/N)$, $f\in C_0(G/N)$
and $r,s\in G$.

Let $\alpha:G\to\Aut A$ be a continuous action of $G$ by automorphisms of
a $C^*$-algebra $A$, and write $\lt$ and $\rt$ for the actions of $G$
on $C_0(G)$ by left and right translation, so that
\[
\lt_s(f)(t)=f(s^{-1}t)\text{\ and\ } \rt_s(f)(t)=f(ts)\text{\ for $f\in C_0(G)$ and $s,t\in G$.}
\]
If $N$ is a closed normal subgroup of $G$, then
there is a natural isomorphism of
$(A\otimes C_0(G/N))\rtimes_{\alpha\otimes\lt}G$
onto $(A\rtimes_{\alpha}G)\rtimes_{\hat\alpha|}G/N$  
(\cite[Lemma~2.3]{EKR-CP}; see also
\cite[Proposition~A.63 and~Theorem~A.64]{BE}).
Representations of both $C^*$-algebras
come from suitably covariant representations $\pi$, $\mu$, and $U$
of $A$, $C_0(G/N)$, and $G$ (respectively) on the same Hilbert space; 
the isomorphism carries $(\pi\otimes\mu)\rtimes U$
to $(\pi\rtimes U)\rtimes \mu$ and for this reason we refer to 
it (and related maps) as the \emph{canonical isomorphism}.

If $A$ and $B$ are $C^*$-algebras, a right-Hilbert $A$--$B$ bimodule is a right Hilbert $B$-module $X$ together with a homomorphism $\phi$ of $A$ into the $C^*$-algebra $\L(X)$ of adjointable operators on $X$; in practice, we suppress $\phi$ and write $a\cdot x$ for $\phi(a)x$. As in \cite{BE}, we view  a right-Hilbert $A$--$B$ bimodule $X$ as a morphism from $A$ to $B$, and say that the diagram
\begin{equation*}
\xymatrix{
A
\ar[rr]^X
\ar[d]_Z
&&
B
\ar[d]^Y
\\
C
\ar[rr]^W
&&
D
}
\end{equation*}
commutes if  $X\otimes_B Y$ and $Z\otimes_C W$ are isomorphic as right-Hilbert $A$--$D$ bimodules. If $\phi:A\to C$ and $\psi:B\to D$ are isomorphisms, then the right-Hilbert $C$--$D$ bimodule $X'$ obtained from $X$ by adjusting the coefficient algebras using $\phi$ and $\psi$ is by definition the bimodule such that
the diagram
\begin{equation*}
\xymatrix{
A
\ar[rr]^X
\ar[d]_\phi^\cong
&&
B
\ar[d]_\cong^\psi
\\
C
\ar[rr]^{X'}
&&
D
}
\end{equation*}
commutes.  (Formally, the left vertical arrow, for example, is the $A$--$C$ bimodule $A$ with $a\cdot b=ab$, $\langle a\,,b\rangle_C=\phi(a^*b)$ and $a\cdot c=a\phi^{-1}(c)$ for $a,b\in A$ and $c\in C$.)
If $B$ is contained in the multiplier algebra $M(A)$ of $A$, 
we denote by $\Res$
the right-Hilbert $B$--$A$ bimodule $A$, where
\[
b\cdot a=ba,\quad a\cdot c=ac,\quad \langle a\,,\, c\rangle_A=a^*c\quad\text{and\ } {}_B\langle a\,,\, c\rangle d=ac^*d
\]
for $a,c,d\in A$ and $b\in B$.

We'll often write ${}_*\langle\cdot\,,\,\cdot\rangle$ and
${}\langle\cdot\,,\,\cdot\rangle_*$ for the left- and right-inner
products, respectively,  in an imprimitivity bimodule, and trust that
it is clear from context in which algebra the values lie.

%======================================================================

\section{Regularity}
\label{reg}

In the coaction context, \emph{regularity} means that the regular
representations are, up to unitary equivalence, precisely those induced
from the trivial quotient group $G/G$: 

\begin{prop}\label{reg-prop}
Let $\delta\colon B\to M(B\otimes C^*(G))$ be a maximal coaction
of a locally compact group $G$ on a $C^*$-algebra $B$. 
Then for each nondegenerate representation $\pi$
of $B$ on a Hilbert space $\H$, 
the representation $\Ind_{G/G}^G(\pi)$ of $B\rtimes_\delta G$
induced using the Mansfield bimodule $Y_{G/G}^G(\delta)$
is unitarily equivalent to the regular representation 
$((\pi\otimes\lambda)\circ\delta)\rtimes(1\otimes M)$ 
of $B\rtimes_\delta G$ on $\H\otimes L^2(G)$.
\end{prop}

Since $\delta$ is a maximal coaction of $G$ on $B$,
by definition of maximality (\cite[Definition~3.1]{EKQ-MC}),
the canonical surjection
\begin{equation}\label{eq-maximal-iso}
(\id\otimes\lambda)\circ\delta\rtimes(1\otimes M)\rtimes(1\otimes\rho)
\colon B\rtimes_\delta G\rtimes_{\hat\delta}G\to
B\otimes\K(L^2(G))
\end{equation}
is an isomorphism.
This makes the $B\otimes\K(L^2(G))-B$ imprimitivity bimodule
$B\otimes L^2(G)$ into 
a $(B\rtimes_\delta G\rtimes_{\hat\delta} G)-B$ 
imprimitivity bimodule which we call the \emph{Katayama bimodule}
(\cite[Definition~4.1]{KQ-MI}),
and which we denote by $K(\delta)$.
By \cite[Corollary~6.2]{KQ-MI}, the imprimitivity bimodules
$K(\delta)$ and $Y_{G/G}^G(\delta)$ are isomorphic, 
so to prove the proposition it suffices to deal with $K(\delta)$.

\begin{proof}
It is straightforward to see that the map $\theta$ determined by
\[
\theta(b\otimes\xi\otimes h)=\pi(b)h\otimes\xi,
\text{\ where\ } b\in B,\ h\in\H,\ \xi\in L^2(G)
\]
extends to a unitary isomorphism of $K(\delta)\otimes_B \H$ 
onto  $\H\otimes L^2(G)$.

Denote by $j_B$ and $j_{C(G)}$ the canonical maps of $B$
and $C_0(G)$ into $M(B\rtimes_\delta G)$. 
To see that~$\theta$ intertwines the induced representation
and the regular representation, 
it suffices to check that
\begin{itemize}
\item[(1)]
$\theta\big(\Ind_{G/G}^G(\pi)(j_B(b))\zeta\big)
=\big(\pi\otimes\lambda(\delta(b))\big)\theta(\zeta)$,
and
\smallskip
\item[(2)] 
$\theta\big(\Ind_{G/G}^G(\pi)(j_{C(G)}(f))\zeta\big)
=(1\otimes M(f))\theta(\zeta)$
\end{itemize}
for $b\in B$, $f\in C_0(G)$, and $\zeta\in K(\delta)\otimes_B \H$; it
further suffices to consider $\zeta$ of the form $a\otimes\eta\otimes h$
for $a\in B$, $\eta\in L^2(G)$ and $h\in \H$.  
Verifying~(2) is straightforward. 
To check~(1), we use nondegeneracy to write
$\eta=\lambda(c)\xi$ for $c\in C^*(G)$ and $\xi\in L^2(G)$; then
$\delta(b)(1\otimes c)\in B\otimes C^*(G)$, and we can approximate it by
a sum $\sum_{j=1}^n b_j\otimes c_j\in B\otimes C^*(G)$.  Now we can do
an approximate calculation:
\begin{align*}
\theta\big(\Ind_{G/G}^G(\pi)(i_B(b))(a\otimes\lambda(c)\xi\otimes h)\big)
&=\theta\big((\id\otimes\lambda(\delta(b))(a\otimes\lambda(c)\xi))\otimes
h\big)\\
&=\theta\big((\id\otimes\lambda(\delta(b)(1\otimes c))(a\otimes
\xi))\otimes h\big)\\
&\sim\sum_{j=1}^n\theta\big((\id\otimes\lambda(b_j\otimes
c_j)(a\otimes\xi))\otimes h\big)\\
&=\sum_{j=1}^n \pi(b_ja)h\otimes\lambda(c_j)\xi\\
&=\pi\otimes\lambda\Big( \sum_{j=1}^n b_j\otimes c_j
\Big)(\pi(a)h\otimes\xi)\\
&\sim\pi\otimes\lambda(\delta(b)(1\otimes c))(\pi(a)h\otimes\xi)\\
&=\pi\otimes\lambda(\delta(b))(\pi(a)h\otimes\lambda(c)\xi)\\
&=\pi\otimes\lambda(\delta(b))\theta(a\otimes\lambda(c)\xi\otimes h);
\end{align*}
since the approximations can be made arbitrarily accurate, 
this implies~(1).
\end{proof}

%======================================================================

\section{A Calculus for Symmetric-Imprimitivity Bimodules}
\label{PhashQ}

The set-up for the 
symmetric imprimitivity theorem of \cite{rae-sit}
is that of commuting free and proper
actions of locally compact groups $K$ and $L$ on the left and right,
respectively, of a locally compact space $P$. In addition, there are
commuting actions $\sigma$ and $\eta$ of $K$ and $L$ on a $C^*$-algebra
$A$.  We sum up this set-up by saying that $({}_KP_L,A,\sigma,\eta)$
is \emph{symmetric imprimitivity data}, and we represent this
schematically with the diagram
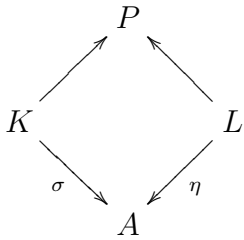
\[
\xymatrix{
& P & \\
K\ar[ru]\ar[dr]_\sigma&& L\ar[lu]\ar[dl]^\eta\\
& A & }
\]
The \emph{induced algebra}
$\Ind_L^P\eta$  consists of all functions $f\in C_b(P, A)$
such that
\[
f(p\cdot t)=\eta^{-1}_t(f(p))
\text{\ for $t\in L$ and $p\in P$  and } 
(pL\mapsto\|f(p)\|)\in C_0(P/L).
\]
Similarly, $\Ind_K^P\sigma$ consists of all functions $f\in C_b(P, A)$ such that
\[
f(s\cdot p)=\sigma_s(f(p))
\text{\ for $s\in K$ and $p\in P$  and  } 
(Kp\mapsto\|f(p)\|)\in C_0(K\backslash P).
\]
$\Ind_L^P\eta$  admits the diagonal action
$\sigma\otimes\lt$ of $K$, and 
$\Ind_K^P\sigma$  admits the diagonal action
$\eta\otimes\rt$ of $L$. 
The symmetric imprimitivity theorem (\cite[Theorem~1.1]{rae-sit})
says that $C_c(P,A)$ can be completed to a
\[
(\Ind_L^P\eta\rtimes_{\sigma\otimes\lt}K)-(\Ind_K^P\sigma\rtimes_{\eta\otimes\rt}L)
\]
imprimitivity bimodule. We denote this bimodule by 
$W({}_KP_L,A,\sigma,\eta)$, or more compactly, by~$W(P)$.

In this section 
we consider two  sets of symmetric imprimitivity
data, $({}_KP_L, A,\sigma,\eta)$ and $({}_LQ_G,A,\xi,\tau)$, which
are compatible in a way that ensures there is an isomorphism
$\Phi$ of $\Ind_K^P\sigma \rtimes_{\eta\otimes\rt}L$ onto
$\Ind_G^Q\tau\rtimes_{\xi\otimes\lt}L$.  
Thus we can form the imprimitivity bimodule $W(P)\otimes_\Phi W(Q)$,
which is by definition the imprimitivity bimodule such that the diagram
\begin{equation*}\label{eq-defn-tensoroverphi}
\xymatrix{
\Ind_L^P\eta\rtimes_{\sigma\otimes\lt}K
\ar[rr]^{W(P)\otimes_\Phi W(Q)}
\ar[d]_{W(P)}
&&
\Ind_L^Q\zeta\rtimes_{\tau\otimes\rt}G
\\
\Ind_K^P\sigma\rtimes_{\eta\otimes\rt}L
\ar[rr]^{\Phi}_\cong
&&
\Ind_G^Q\tau\rtimes_{\zeta\otimes\lt}L
\ar[u]_{W(Q)}
}
\end{equation*}
of imprimitivity bimodules commutes.
Theorem~\ref{thm-PhashQ} will show that
$W(P)\otimes_\Phi W(Q)$ can be replaced with an
imprimitivity bimodule based on a single set of symmetric imprimitivity
data, thus giving an easy way of calculating, at the level of spaces,
the isomorphism class of the balanced tensor product.

Suppose 
$\phi:K\under P\to Q/G$ is a homeomorphism 
which is $L$-equivariant in the sense that
\begin{equation}\label{L-eqvt}
\phi(K\cdot p\cdot t)=t^{-1}\cdot\phi(K\cdot p)
\quad\text{for all $t\in L$ and $p\in P$,} 
\end{equation}
and let  $P\times_{\phi}Q:=\{(p,q)\in P\times Q:\phi(K\cdot p)=q\cdot
G\}$
be the fibred product. 
We define 
\begin{equation}
  \label{eq:4}
  \pspq:=(P\times_{\phi}Q)/L,
\end{equation}
where the action of $L$ on $P\times_{\phi}Q$ is via the diagonal action
$(p,q)\cdot t:=(p\cdot t,t^{-1}\cdot q)$.
We will use $[p,q]$ to denote the class of $(p,q)$ in $\pspq$;
we will write $\psq$ for $\pspq$ when there is no risk of
confusion. 

\begin{thm}
\label{thm-PhashQ}
Suppose $K$, $L$ and $G$ are locally compact groups, 
and suppose that $({}_KP_L,A,\sigma,\eta)$
and $({}_LQ_G,A,\zeta,\tau)$ are symmetric imprimitivity data.
In addition, suppose there is an $L$-equivariant homeomorphism
$\phi:K\under P\to Q/G$ as at~\eqref{L-eqvt}, 
and that there are continuous maps
\begin{equation*}
\sigmat:P\to\Aut A\quad\text{and}\quad \taut:Q\to\Aut A
\end{equation*}
such that, for $p\in P$, $q\in Q$, $k\in K$, $m\in G$ and $t\in L$,
\begin{gather}
\label{eq:5}
\sigmat_{k\cdot p\cdot t}=\sigma_{k} \sigmat_{p} \zeta_{t},
\\
\label{eq:6} \taut_{t\cdot q\cdot m}= \eta_{t}  \taut_{q}
\tau_{m}\text{, and} \\
\text{$\zeta$, $\sigma$ and $\sigmat$ commute with $\eta$, $\taut$ and
$\tau$.} \label{eq:8}
\end{gather}
Then $\pspq$, as defined at~\eqref{eq:4}, admits 
commuting free and proper  actions of  $K$
and $G$, and there are isomorphisms
\begin{gather}
\label{eq:7}
\Phi:\indkp\sigma\rtimes_{\eta\otimes\rt}L\to
\indgq\tau\rtimes_{\zeta\otimes\lt}L , \\
\Phi_{\sigma}:\indkpsq\sigma\rtimes_{\tau\otimes\rt}G\to
\indlq\zeta\rtimes_{\tau\otimes \rt}G \text{, and}\label{eq:10}\\
\Phi_{\tau}:\indgpsq\tau\rtimes_{\sigma\otimes\lt}K \to
\indlp\eta\rtimes_{\sigma\otimes\lt}K,\label{eq:11}
\end{gather}
such that the diagram
\begin{equation}\label{eq-PhasQthm}
\xymatrix{
\Ind_G^{P\#Q}\tau\rtimes_{\sigma\otimes\lt}K
\ar[d]_{\Phi_\tau}
^\cong
\ar[rr]^{W(P\#Q)}
&& \Ind_K^{P\#Q}\sigma\rtimes_{\tau\otimes\rt}G
\ar[d]^{\Phi_\sigma}_\cong
\\
\Ind_L^P\eta\rtimes_{\sigma\otimes\lt}K
\ar[rr]^{W(P)\otimes_\Phi W(Q)}
&& \Ind_L^Q\zeta\rtimes_{\tau\otimes\rt}G}
\end{equation}
of imprimitivity bimodules commutes.
\end{thm}

In the proof of Theorem~\ref{thm-PhashQ} we will use the 
following lemma
to establish the isomorphisms \eqref{eq:7}--\eqref{eq:11}.
In part~(2) of the lemma, $\phi_*$ denotes the natural isomorphism of 
$C_0(Q/G,A)$ onto $C_0(\kup,A)$ induced by $\phi$.

\begin{lem}\label{lem-one-and-two}
Assume the hypotheses of Theorem~\ref{thm-PhashQ}.
\begin{enumerate}
\item There are isomorphisms $\phi_{\sigma}:\indkpsq\sigma \to
\indlq\zeta$ and $\phi_{\tau}:\indgpsq\tau\to\indlp\eta$ given by
\begin{align*}
\phi_{\sigma}(f)(q)=\sigmat^{-1}_{p}\bigl(f\bigl([p,q]\bigr)\bigr)
\quad\text{and}\quad
\phi_{\tau}(f)(p)=\taut_{q}\bigl(f\bigl([p,q]\bigr)\bigr),
\end{align*}
where $p\in P$ and $q\in Q$ are such that $\phi(K\cdot p)=q\cdot G$.
These  isomorphisms are equivariant and hence induce isomorphisms
$\Phi_{\sigma}:=\phi_{\sigma}\rtimes G$ and $\Phi_{\tau}:=
\phi_{\tau}\rtimes K$ of the crossed products.
\item The maps defined by
\begin{equation*}
\psi_{\sigma}(f)(K\cdot p):= \sigmat^{-1}_{p}\bigl(f(p)\bigr)
\quad\text{and} \quad \psi_{\tau}(f)(q\cdot G):= \taut_{q} \bigl(f(q)\bigr)
\end{equation*}
give isomorphisms $\psi_{\sigma}:\indkp\sigma\to C_{0}(\kup,A)$ and
$\psi_{\tau} : \indgq\tau\to C_{0}(\qmm,A)$.  Furthermore, the
composition
\begin{equation*}
\xymatrix{\indkp\sigma\ar[r]^-{\psi_{\sigma}}
\ar@/_2pc/[rrr]_{T}
&C_{0}(\kup,A)
\ar[r]^-{\phi_{*}^{-1}} & C_{0}(\qmm,A)\ar[r]^-{\psi_{\tau}^{-1}}&\indgq\tau}
\end{equation*}
is given by
\begin{equation}\label{Tfq}
T(f)(q)=\taut_{q}^{-1}
\sigmat_{p}^{-1}\bigl(f(p)\bigr),
\end{equation}
where $p\in P$ is such that $\phi(K\cdot p)=q\cdot G$.
$T$ is equivariant, and hence $\Phi:=T\rtimes L$ is an
isomorphism of $\indkp\sigma\rtimes_{\eta\otimes\rt}L$ onto
$\indgq\tau\rtimes_{\zeta\otimes \lt}L$.
\end{enumerate}
\end{lem}

\begin{proof}
(1) The first step is to verify that $\phi_{\sigma}$ is well-defined. Let $f\in\indkpsq \sigma$.
If $\phi(K\cdot p)=q\cdot G$, then for any $k\in K$,
\[
\sigmat_{k\cdot p}\bigl(f\bigl([k\cdot p,q]\bigr)\bigr)
=(\sigma_{k}
\sigmat_{p})^{-1}\bigl(\sigma_{k}\bigl(f\bigl([p,q]\bigr)\bigr)\bigr)
= \sigmat_{p}^{-1}\bigl(f\bigl([p,q]\bigr)\bigr).
\]
It follows that $\phi_{\sigma}(f)$ is a well-defined function on $Q$.  On
the other hand, if $\phi(K\cdot p)=q\cdot G$, then, for all $t\in L$,  $\phi(K\cdot p\cdot
t^{-1}) = t\cdot q\cdot G$ and
\begin{align*}
\phi_{\sigma}(f)(t\cdot q) &= \sigmat_{p\cdot
t^{-1}}^{-1}\bigl(f\bigl([p\cdot t^{-1},t\cdot q]\bigr)\bigr)
=\sigmat_{p\cdot t^{-1}}^{-1}\bigl(f\bigl([p,q]\bigr)\bigr) \\
&= \zeta_{t}  \sigmat_{p}^{-1}\bigl(f\bigl([p,q]\bigr)\bigr)
= \zeta_{t}\bigl(\phi_{\sigma}(f)(q)\bigr).
\end{align*}
Therefore, to see that $\phi_{\sigma}(f)$ is in $\indlq\zeta$, we only
have to check that $\phi_{\sigma}(f)$ is continuous and that $L\cdot
q\mapsto \|f(q)\|$ vanishes at infinity.

To establish continuity, it suffices to show that, given any net $q_{\alpha}\to q$
we can find a subnet such that, after we pass to the subnet and
relabel, we have $\phi_{\sigma}(f)(q_{\alpha})\to
\phi_{\sigma}(f)(q)$.  Choose $p_{\alpha}$ such that $\phi(K\cdot
p_{\alpha})=q_{\alpha}\cdot G$.  Since $\phi$ is a homeomorphism,
there is a $p$ such that
\begin{equation*}
K\cdot p_{\alpha}\to K\cdot p=\phi^{-1}(q\cdot G).
\end{equation*}
Since the orbit map is open, we can pass to a subnet, relabel, and
assume that there are $k_{\alpha}\in K$ such that
\begin{equation*}
k_{\alpha}\cdot p_{\alpha}\to p.
\end{equation*}
Of course, $\phi(Kk_{\alpha}\cdot p_{\alpha})=q_\alpha\cdot G$, and
\begin{align*}
  \phi_{\sigma}(f)(q_{\alpha})&= \sigmat_{k_{\alpha}\cdot
  p_{\alpha}}^{-1} \bigl(f\bigl([k_{\alpha}\cdot p_{\alpha},
  q_{\alpha}]\bigr)\bigr)
\to \sigmat_{p}^{-1}\bigl(f\bigl([p,q]\bigr)\bigr)
=\phi_{\sigma}(f)(q)
\end{align*}
because $f$ and $\sigmat$ are continuous.
Thus, $\phi_{\sigma}(f)$ is continuous.

To see that $\phi_{\sigma}(f)$ vanishes at infinity, it suffices to
show that if $\{q_{\alpha}\}$ is a net in $Q$ such that
\begin{equation*}
  \|\phi_{\sigma}(f)(q_{\alpha})\|\ge\epsilon>0,
\end{equation*}
then $\{q_{\alpha}\}$ has a convergent subnet.  Let $p_{\alpha}$ be
such that $\phi(K\cdot p_{\alpha})=q_{\alpha}\cdot G$.  Since
$\sigmat_{p_{\alpha}}^{-1}$ is isometric, we must have
\begin{equation*}
\|f\bigl([p_{\alpha},q_{\alpha}]\bigr)\|\ge\epsilon\quad\text{for all
$\alpha$.}
\end{equation*}
Then, since $K\cdot[p,q]\mapsto \|f\bigl([p,q]\bigr)\|$ vanishes at
infinity, we can pass to a subnet, relabel, and assume that there is a
$[p,q]\in\psq$ such that
\begin{equation*}
K\cdot[p_{\alpha},q_{\alpha}]\to K\cdot[p,q].
\end{equation*}
Since orbit maps are open, we can pass to another subnet, relabel, and
find $k_{\alpha}\in K$ such that
\begin{equation*}
k_{\alpha}\cdot[p_{\alpha},q_{\alpha}]=[k_{\alpha}\cdot p_{\alpha},
q_{\alpha}] \to [p,q].
\end{equation*}
Similarly, after passing to another subnet and relabeling, there are
$t_{\alpha}\in L$ such that
\begin{equation*}
(k_{\alpha}\cdot p_{\alpha}\cdot t_{\alpha},t_{\alpha}^{-1}\cdot
q_{\alpha}) \to (p,q).
\end{equation*}
In particular, $K\cdot q_{\alpha}\to K\cdot q$,  and  hence $\phi_\sigma(f)\in\indlq\zeta$.
Since the operations are pointwise, $\phi_{\sigma}$ is a homomorphism of $\indkpsq\sigma$ into
$\indlq\zeta$; it is an isomorphism since similar considerations show
that
\begin{equation*}
\phi_{\sigma}^{-1}(g)([p,q])=\sigmat_{p}\bigl(g(q)\bigr)
\end{equation*}
is an inverse.

Furthermore, if $m\in G$, and if $\phi(K\cdot p)=q\cdot G$, then
$\phi(K\cdot p)=(q\cdot m)\cdot G$, and, since $\tau$ and $\sigmat$ commute,
\begin{align*}
  \sigmat_{p}^{-1}\bigl(\tau_{m}\otimes\rt_{m}(f)\bigl([p,q]\bigr)\bigr)
&=\sigmat_{p}^{-1} \bigl(\tau_{m}\bigl(f\bigl([p,qm]\bigr)\bigr)\bigr)
  \\
&=\tau_{m}\bigl(\sigmat_{p}^{-1}\bigl(f\bigl([p,qm]\bigr)\bigr)\bigr)
  \\
&= \tau_{m}\bigl(\phi_{\sigma}(f)(q\cdot m)\bigr) \\
&= (\tau_{m}\otimes\rt_{m})\phi_{\sigma}(f)(q).
\end{align*}
Thus $\phi_{\sigma}$ is equivariant, and therefore
gives an isomorphism
$\Phi_{\sigma}=\phi_{\sigma}\rtimes G$.

The statements for $\phi_{\tau}$ and $\Phi_{\tau}$ are proved
similarly.

(2)
It is easy to check that $\psi_{\sigma}$ and $\psi_{\tau}$
are well-defined homomorphisms which are isomorphisms
by computing their inverses directly (for example,
$\psi_\tau^{-1}(g)(q)=\widetilde{\tau_q}^{-1}(g(q\cdot G))$),
and it is then straightforward to verify~\eqref{Tfq}.  
Further, if $\phi(K\cdot
p)=q\cdot G$ and if $t\in L$, then on the one hand
\begin{equation}\label{eq:x3}
T\bigl(\eta_{t}\otimes\rt_{t}(f)\bigr)(q)= \taut_{q}^{-1}
\sigmat_{p}^{-1} \bigl(\eta_{t}\otimes \rt_{t}(f)(p)\bigr) 
= \taut_{q}^{-1}
\sigmat_{p}^{-1}\bigl(\eta_{t}\bigl(f(p\cdot t)\bigr)\bigr).
\end{equation}
On the other hand, we also have $\phi(K\cdot p\cdot t)=t^{-1}\cdot
q\cdot G$, and
\begin{align*}
\zeta_{t}\otimes\lt_{t}\bigl(T(f)\bigr)(q) &=
\zeta_{t}\bigl(T(f)(t^{-1} \cdot q)\bigr)\\
&= \zeta_{t}\bigl(\taut_{t^{-1}\cdot q}^{-1}  \sigmat_{p\cdot
t}^{-1} \bigl(f(p\cdot t)\bigr)\bigr)\\
\intertext{which, since $\taut_{t^{-1}\cdot q}=\eta_{t}^{-1}
\taut_{q}$ and $\sigma_{p\cdot t}=\sigmat_{p}  \zeta_{t}$, is}
&= \zeta_{t}\bigl(\taut_{q}^{-1}  \eta_{t}  \zeta_{t}^{-1}
\sigmat_{p} \bigl(f(p\cdot t)\bigr)\bigr),
\end{align*}
and this coincides with \eqref{eq:x3} because $\zeta$ commutes with $\taut$
and $\eta$, and $\eta$ commutes with $\sigmat$.
Thus, $T$ is equivariant and the result follows.
\end{proof}

\begin{proof}[Proof of Theorem~\ref{thm-PhashQ}]
Let $\Phi_\sigma$, $\Phi_\tau$ and $\Phi$ be as in Lemma~\ref{lem-one-and-two}.  For fixed  $x\in \cc(P,A)$, $y\in\cc(Q,A)$ and $(p,q)\in P\times Q$ set
\begin{equation}\label{eq:14}
f(p,q):=\int_{L} \taut_{r^{-1}\cdot q}^{-1}\bigl(x(p\cdot r)\bigr)
\sigmat_{p\cdot r} \bigl(y(r^{-1}\cdot q)\bigr)\,dr.
\end{equation}
Straightforward computation using the left-invariance of Haar
measure shows that $f(p,q)$ depends only on the class $[p,q]$ of $(p,q)\in\psq$. Since the actions of $L$ on $P$ and $Q$ are free and proper,  $f(p,q)<\infty$ and $[p,q]\mapsto f(p,q)$ is continuous with compact support.
Thus we can define $\Omega:\cc(P,A)\atensor\cc(Q,A)\to
\cc(\psq,A)$ by letting
\[\Omega(x\otimes y)([p,q])=f(p,q).\]
(That $\Omega$ is well-defined on the balanced tensor product
will follow from the same calculation that shows $\Omega$ is isometric for the right inner products; see below.)

To see that \eqref{eq-PhasQthm} commutes, we will show that the triple
$(\Phi_\tau^{-1},\Omega,\Phi_\sigma^{-1})$ extends to an imprimitivity
bimodule isomorphism of $W(P)\otimes_\Phi W(Q)$ onto $W(P\# Q)$.  In
particular, we will show that $(\Phi_\tau^{-1},\Omega,\Phi_\sigma^{-1})$
preserves the right inner products and both the left
and right actions.  Then the range of $\Omega$ will be a closed
sub-bimodule of $\ZZ(\psq)$ on which the right inner product is full.
It will then follow from  the Rieffel correspondence (see, for example,
\cite[Proposition~3.24]{tfb}) that $\Omega$ is surjective.  This will
imply that $(\Phi_\tau^{-1},\Omega,\Phi_\sigma^{-1})$  must also preserve
the left inner product and hence will be the desired isomorphism.

Let $x,w\in\cc(P,A)\subseteq W(P)$ and
$y,z\in\cc(Q,A)\subseteq W(Q)$ and let $\langle\langle\cdot\,,\,\cdot\rangle\rangle_*$
be the right inner product on $\ZZ(P)\otimes_{\Phi}\ZZ(Q)$.  We will show that
\[\langle\langle x\otimes y\,,\,w\otimes z\cdot\rangle\rangle_*=\Phi_\sigma\big( \langle\Omega(x\otimes y)\,,\,\Omega(w\otimes z)\rangle_* \big).\]
The inner product $\langle\langle x\otimes y\,,\,w\otimes z\cdot\rangle\rangle_*$ takes values in $\cc(G,\indkq\zeta)\subseteq \indlq\zeta\rtimes G$ which we view
as  functions on $G\times Q$.  Thus
\begin{align*}
 \big\langle\langle x&\otimes y\,,\,w\otimes z\rangle\big\rangle_*(m,q)
=\big\langle \Phi\bigl(\langle w\,,\,x\rangle_*\bigr)\cdot y\,,\,z\big\rangle_* (m,q) \\
&=
\Delta_G(m)^{1/2}\int_{L} \zeta_{t} \bigl(\bigl(\Phi(\langle w\,,\, x\rangle_*)\cdot y\bigr)(t^{-1}
\cdot q)^{*}\tau_{m}\bigl(z(t^{-1}\cdot q\cdot m)\bigr)\bigr) \,dt \\
&=
 \Delta_G(m)^{1/2}\int_{L}\zeta_{t}\Bigl(\int_{L}\Phi\bigl(\langle w\,,\,x\rangle_*\bigr) (r,t^{-1}\cdot q)
\zeta_{r} \bigl(y(r^{-1}t^{-1}\cdot q)\bigr)\Delta_L(r)^{1/2} \,dr\Bigr)^{*}
\\
&\hskip 10 cm
\zeta_{t}  \tau_{m} \bigl(z(t^{-1}\cdot q\cdot m)\bigr) \,dt\\
\intertext{which, if $\phi(K\cdot p)=q\cdot G$, is}
&=
\Delta_G(m)^{-1/2}\int_{L}\int_{L} \zeta_{t}\Bigl( \taut_{t^{-1}\cdot q}^{-1}
\sigmat_{p\cdot t}^{-1} \bigl(\langle w\,,\,x\rangle_*(r,p\cdot t)\bigr)
\zeta_{r}\bigl( y(r^{-1}t^{-1}\cdot q)\bigr) \Bigr)^{*}\Delta_L(r)^{1/2}\,dr \\
&\hskip11 cm
{\zeta_{t}  \tau_{m}\bigl(z(t^{-1}\cdot q\cdot m)\bigr) }\,dt\\
&=
\Delta_G(m)^{-1/2}\int_{L}\int_{L}
\zeta_{tr} \bigl(y(r^{-1}t^{-1}\cdot q)\bigr)^{*}\\
&\hskip1.6 cm\zeta_{t} \taut_{t^{-1}\cdot q}^{-1}\sigmat_{p\cdot t}^{-1}
\Bigl( \int_{K}\sigma_{s}\bigl(w(s^{-1}\cdot p\cdot t)^{*} \eta_{r}
\bigl( x(s^{-1}\cdot p\cdot tr)\bigr)\bigr)\,ds\Bigr)^{*}
\zeta_{t}
\tau_{m} \bigl( z(t^{-1}\cdot q\cdot m)\bigr) \,dr\,dt \\
\intertext{which, since $\taut_{t^{-1}\cdot
q}=\eta_{t}^{-1} \taut_{q}$ and $\sigmat_{p\cdot
t}=\sigmat_{p} \zeta_{t}$, and since $\zeta$ commutes with both $\taut$
and $\eta$, and $\eta$ commutes with $\sigmat$ (see \eqref{eq:5}--\eqref{eq:8}), is}
&=
\Delta_G(m)^{-1/2}\int_{L}\int_{L}\int_{K} \zeta_{tr}\bigl(y(r^{-1}t^{-1}\cdot
q)^{*}\bigr)
\\
&\hskip2 cm\taut_{q}^{-1} \sigmat_{p}^{-1} \Bigl( \sigma_{s} \eta_{tr}
\bigl( x(s^{-1}\cdot p \cdot tr)^{*}\bigr) \eta_{t}  \sigma_{s}
\bigl( w(s^{-1}\cdot p\cdot t)\bigr)\Bigr)
\zeta_{t}\tau_{m} \bigl(z(t^{-1}\cdot q\cdot m)\bigr)
\,ds\,dr\,dt \\
\intertext{which, replacing $r$ by $t^{-1}r$ and using
\eqref{eq:5}--\eqref{eq:8} again, is}
&=
 \Delta_G(m)^{-1/2}\int_{K}\Bigl(\int_{L} \zeta_{r}\bigl(y(r^{-1}\cdot q)^{*}\bigr)
\taut_{r^{-1}\cdot q}^{-1} \sigmat_{s^{-1}\cdot p}^{-1}\bigl(
x(s^{-1}\cdot p\cdot r)^{*}\bigr)\Bigr) \\
&\hskip 4 cm
\Bigl( \int_{L} \taut_{t^{-1}\cdot q}^{-1}\sigmat_{s^{-1}\cdot
  p}^{-1} \bigl(w(s^{-1}\cdot p \cdot t)\bigr) \zeta_{t}\tau_{m}
\bigl( z(t^{-1}\cdot q\cdot m)\bigr)\,dt \Bigr)\,ds\\
\intertext{which, using \eqref{eq:5}, is}
&=
\Delta_G(m)^{-1/2}\int_{K} \sigmat_{s^{-1}\cdot p}^{-1} \Bigl(
\Bigl(\int_{L}\taut_{r^{-1}\cdot q}^{-1}\bigl(x(s^{-1}\cdot p \cdot
r)\bigr)
 \sigmat_{s^{-1}\cdot p\cdot r} \bigl(y(r^{-1}\cdot
q)\bigr)\,dr\Bigr)^{*} \\
&\hskip4cm
\Bigl(\int_{L}\taut_{t^{-1}\cdot q}^{-1}\bigl(w(s^{-1}\cdot p\cdot
t)\bigr) \sigmat_{s^{-1}\cdot p\cdot
 t}\tau_{m}\bigl(z(t^{-1}\cdot q\cdot m)\bigr)\,dt\Bigr)
\Bigr)\,ds \\
\intertext{which, using \eqref{eq:8} and the definition of $\Omega$, is}
&=
\Delta_G(m)^{-1/2}\int_{K}\sigmat_{s^{-1}\cdot p}^{-1} \Bigl(
\Omega(x\otimes y)\bigl([s^{-1}\cdot p,  q]\bigr)^{*} \\
&\hskip7em
\tau_{m} \Bigl(\int_{L} \taut_{t^{-1}\cdot q\cdot m}^{-1} \bigl(
w(s^{-1} \cdot p \cdot t)\bigr)\sigmat_{s^{-1}\cdot p\cdot t} \bigl(
z(t^{-1} \cdot q\cdot m)\bigr)\,dt \Bigr)
\Bigr)\,ds \\
&=
\Delta_G(m)^{-1/2}\sigmat_{p}^{-1} \Bigl(
\int_{K} \sigma_{s} \bigl(\Omega(x\otimes y)\bigl([s^{-1}\cdot
p,q]\bigr)^{*} \tau_{m}\bigl(\Omega(w\otimes z)\bigl([s^{-1}\cdot
p,q\cdot m]\bigr)\bigr)\,ds\bigr)
\Bigr)\\
&=
\sigmat_{p}^{-1} \bigl(\langle\Omega(x\otimes y)\,,\,\Omega(w\otimes z)\rangle_*
(m,[p,q]\bigr)\bigr) \\
&= \Phi_{\sigma}\bigl(\langle\Omega(x\otimes y)\,,\,\Omega(w\otimes
z)\rangle_*\bigr) (m,q).
\end{align*}
Thus $(\Phi_\tau^{-1},\Omega,\Phi_\sigma^{-1})$ intertwines the right inner products.

If $b\in\cc(G,\indkpsq\sigma)\subseteq\Ind_K^{P\#Q}\sigma\rtimes_{\tau\otimes\rt}G$ is viewed as a function on $G\times (\psq)$,
then
\begin{align}
\Omega(x&\otimes y)\cdot b \bigl([p,q]\bigr) =
\int_{G} \tau_{m}\bigl(\Omega(x\otimes y)\bigl([p,q\cdot m]\bigr)\bigr)
b(m^{-1},[p,q\cdot m])\bigr)\bigr)\Delta_G(m)^{-1/2}\,dm \notag\\
&=
\int_{G} \tau_{m} \Bigl(
\int_{L}\taut_{r^{-1}\cdot q\cdot m}^{-1}\bigl(x(p\cdot r)\bigr)
\sigmat_{p\cdot r} \bigl(y(r^{-1}\cdot q\cdot m)\bigr) \,dr
\Bigr) \notag
%\\
%&\hskip10em
b(m^{-1},[p,q\cdot m])\Delta_G(m)^{-1/2}\,dm \notag\\
&=
\int_{G}\int_{L} \taut_{r^{-1}\cdot q}^{-1}\bigl(x(p\cdot r)\bigr)
\tau_{m}   \sigmat_{p\cdot r}\bigl(y(r^{-1}\cdot q\cdot m)\bigr)
\label{eq:12}
%\\
%&\hskip10em
b(m^{-1},[p,q\cdot m])\,dr\Delta_G(m)^{-1/2}\,dm.\notag
\end{align}
On the other hand,
\begin{align*}
\Omega\bigl((&x\otimes y)\cdot \Phi_{\sigma}(b))\bigl([p,q]\bigr)
=
\Omega(x\otimes
(y\cdot \Phi_{\sigma}(b))\bigl([p,q]\bigr) \\
&=
\int_{L} \taut_{r^{-1}\cdot q}^{-1} \bigl(x(p\cdot r)\bigr)
\sigmat_{p\cdot r} \bigl(y\cdot \Phi_{\sigma}(b)(r^{-1}\cdot
q)\bigr) \,dr \\
&=
\int_{L} \taut_{r^{-1}\cdot q}^{-1} \bigl(x(p\cdot r)\bigr)
\sigmat_{p\cdot r} \Bigl(\int_G\tau_m\bigl(y(r^{-1}\cdot q\cdot m)\Phi_{\sigma}(b)(m^{-1},r^{-1}\cdot
q\cdot m)\bigr)\Delta_G(m)^{-1/2}\, dm\Bigr) \,dr \\
&=\int_{L}\int_{G} \taut_{r^{-1}\cdot q}^{-1}\bigl(x(p\cdot r)\bigr)
\sigmat_{p\cdot r}\tau_{m} \bigl(y(r^{-1}\cdot q\cdot m)\bigr)
b(m^{-1},[p,q\cdot m]) \Delta_G(m)^{-1/2}\,dm\,dr
\end{align*}
where we have used that $\phi(K\cdot p\cdot r)=r^{-1}\cdot q\cdot G$
implies
\[\Phi_{\sigma}(b)(m^{-1},r^{-1}\cdot q\cdot m) =
\sigmat_{p\cdot r}^{-1}\bigl(b(m^{-1},[p\cdot r,r^{-1}\cdot q\cdot
m])\bigr) =\sigmat_{p\cdot r}^{-1}\bigl(b(m^{-1},[p,q\cdot
m])\bigr).\]
Since $\tau$ and $\sigmat$ commute,
an application of Fubini's Theorem gives $\Omega(x\otimes y)\cdot b =\Omega\bigl((x\otimes y)\cdot \Phi_{\sigma}(b)\bigr)$.

For the left action, let
$c\in\cc(K, \indgpsq\tau)\subseteq \indgpsq\alpha\rtimes_{\sigma\otimes\lt} K$.  Then,
viewing $c$ as a function on $K\times (\psq)$, we have
\begin{align}
c\cdot\Omega(x\otimes y)&\bigl([p,q]\bigr)
=
\int_{K} c(t,[p,q])\sigma_{t}\bigl(\Omega(x\otimes
y)\bigl([t^{-1}\cdot p,q]\bigr)\bigr) \Delta_K(t)^{1/2}\,dt\notag \\
&=
\int_{K}
c(t,[p,q])\sigma_{t}\Bigl(
\int_{L}\taut_{r^{-1}\cdot q}^{-1}\bigl(x(t^{-1}\cdot p\cdot r)\bigr)
\sigmat_{t^{-1}\cdot p\cdot r} \bigl(y(r^{-1}\cdot q)\bigr)\,dr
\Bigr)\Delta_K(t)^{1/2}\,dt \notag\\
&=
\int_{K}\int_{L} c(t,[p,q])\sigma_{t}  \taut_{r^{-1}\cdot q}^{-1}
\bigl( x(t^{-1}\cdot p \cdot r)\bigr) \sigmat_{p\cdot r}
\bigl(y(r^{-1}\cdot q )\bigr) \Delta_K(t)^{1/2}\,dr\,dt.
\label{eq:15}
\end{align}
On the other hand,
\begin{align*}
\Omega(\Phi_\tau(c)\cdot (x&\otimes y))\bigl([p,q]\bigr) =
\Omega(\Phi_{\tau}(c)\cdot x\otimes y)\bigl([p,q]\bigr) \\
&=
\int_{L} \taut^{-1}_{r^{-1}\cdot q} \bigl(\Phi_{\tau}(c)\cdot x(p\cdot
r)\bigr) \sigmat_{p\cdot r}\bigl(y(r^{-1}\cdot q)\bigr) \,dr \\
&=
\int_{L} \taut_{r^{-1}\cdot q}^{-1} \Bigl(
\int_{K} \Phi_{\tau}(c)(t,p\cdot r)\sigma_{t}\bigl(x(t^{-1}\cdot
p\cdot r)\bigr) \Delta_K(t)^{1/2}\,dt
\Bigr)
\sigmat_{p\cdot r} \bigl(y(r^{-1}\cdot q)\bigr) \,dr  \\
&=
\int_{L} \taut_{r^{-1}\cdot q}^{-1} \Bigl(
\int_{K} \taut_{r^{-1}\cdot q}\bigl(c(t,[p,q])\bigr)
\sigma_{t}\bigl(x(t^{-1}\cdot p\cdot r)\bigr)\Delta_K(t)^{1/2}\,dt
\Bigr)
\sigmat_{p\cdot r} \bigl(y(r^{-1}\cdot q)\bigr) \,dr \\
&=
\int_{L}\int_{G} c(t,[p,q])\taut_{r^{-1}\cdot q}^{-1}  \sigma_{t}
\bigl(x(t^{-1}\cdot p\cdot r)\bigr) \sigmat_{p\cdot
  r}\bigl(y(r^{-1}\cdot q)\bigr)\Delta_K(t)^{1/2} \,dt\,dr,
\end{align*}
which coincides with \eqref{eq:15}.  This
completes the proof.
\end{proof}
%======================================================================

\section{Imprimitivity Bimodule Isomorphisms}
\label{previous}

In this section, we show that for a dual coaction $\hat\alpha$,
the Mansfield bimodule
$Y_{G/H}^{G/N}(\hat\alpha|)$ appearing in Theorem~\ref{main-thm}
can be replaced by a symmetric imprimitivity bimodule.
More precisely, we show how this result (Theorem~\ref{thm-YZ})
follows from a certain bimodule factorization result
(Theorem~\ref{thm-KXZ}). 
Preparation for the proof of Theorem~\ref{thm-KXZ} takes up the rest of
this section; the proof itself occupies Section~\ref{next}. 

\begin{thm}\label{thm-YZ}
Suppose $\alpha$ is a continuous action of a locally compact group $G$ 
by automorphisms of a $C^*$-algebra $A$, and suppose $N$ and $H$ are closed
normal subgroups of $G$ with $N\subseteq H$.  
Let $\epsilon$ be the maximal coaction $\hat\alpha|_{G/N}$ of $G/N$
on $A\rtimes_\alpha G$.
Then the diagram
\begin{equation}\label{YZ}
\xymatrix{
A\rtimes_{\alpha}G\rtimes_{\epsilon}(G/N)\rtimes_{\hat\epsilon|}(H/N)
\ar[rr]^-{Y_{G/H}^{G/N}(\epsilon)}
\ar[d]^\cong
&& 
A\rtimes_{\alpha}G\rtimes_{\epsilon|}(G/H)
\ar[d]_\cong
\\
(A\otimes C_0(G/N))\rtimes_{\alpha\otimes\lt}G\rtimes_{\beta|}(H/N)
\ar[rr]^-{Z_{G/H}^{G/N}(\alpha)}
&& 
(A\otimes C_0(G/H))\rtimes_{\alpha\otimes\lt}G
}
\end{equation}
of imprimitivity bimodules commutes, where the vertical arrows are 
the canonical isomorphisms.
\end{thm}

Here $Z_{G/H}^{G/N}(\alpha)$ is the symmetric-imprimitivity
bimodule constructed in \cite[Proposition~3.3]{aHKRW-EP};
we will review its construction in Section~\ref{Q-sec}.
The action 
\begin{equation}\label{new-beta}
\beta\colon G/N\to 
\Aut((A\otimes C_0(G/N))\rtimes_{\alpha\otimes\lt}G).
\end{equation}
is induced by the action
$\id\otimes\rt$ of $G/N$ on $A\otimes C_0(G/N)$,
which commutes with the action $\alpha\otimes\lt$ of $G$.
It corresponds to the dual 
of the coaction $\epsilon=\hat\alpha|_{G/N}$
under the canonical isomorphism 
of $A\rtimes_{\alpha} G\rtimes_{\hat\alpha|}(G/N)$
with $(A\otimes C_0(G/N))\rtimes_{\alpha\otimes\lt}G$.

Dual coactions and their restrictions are maximal by 
\cite[Proposition~3.4]{EKQ-MC} and \cite[Corollary~7.2]{KQ-MI}, 
so $\epsilon$ is a maximal coaction of $G/N$ on~$A\rtimes_\alpha G$.
Thus, 
the Katayama bimodule $K(\epsilon)$ (see the discussion
following~\eqref{eq-maximal-iso})
is an 
$(A\rtimes_\alpha G)\rtimes_\epsilon G/N\rtimes_{\hat\epsilon} G/N-
A\rtimes_\alpha G$ 
imprimitivity bimodule.
By \cite[Proposition~4.2]{KQ-MI}, $K(\epsilon)$  comes equipped with a
$\hat{\hat\epsilon}$--$\epsilon$ compatible coaction $\epsilon_K$,
and we can further restrict these coactions to $G/H$ 
and take crossed products (see, for example, \cite[\S3.1.2]{BE}).
%to get a 
%$B\rtimes_{\epsilon} G/N
%\rtimes_{\hat\epsilon}(G/N)
%\rtimes_{\hat{\hat\epsilon}|}(G/H)
%-B\rtimes_{\epsilon|}(G/H)$ 
%imprimitivity bimodule $K(\epsilon)\rtimes_{\epsilon_K|}(G/H)$.

The following factorization theorem for $Z_{G/H}^{G/N}(\alpha)$
generalises \cite[Proposition~6.3]{KQ-MI},
which is the one-subgroup version.

\begin{thm}\label{thm-KXZ}
Under the hypotheses of Theorem~\ref{thm-YZ}, 
%Suppose $\alpha$ is a continuous action of a locally compact group 
%$G$ by automorphisms of a $C^*$-algebra $A$, 
%and suppose $N$ and $H$ are closed normal subgroups of $G$
%with $N\subseteq H$.  
%Let $\epsilon$ be the maximal
%coaction $\hat\alpha|_{G/N}$ of $G/N$ on
%$A\rtimes_\alpha G$.
%Then 
the diagram 
\begin{equation}
\label{eq-KXZ}
\xymatrix{
((A\otimes C_0(G/N))\rtimes_{\alpha\otimes\lt}G\rtimes_{\beta|}(H/N)
\ar[r]^-{Z_{G/H}^{G/N}(\alpha)}
&
(A\otimes C_0(G/H))\rtimes_{\alpha\otimes\lt}G
\\
A\rtimes_\alpha G\rtimes_{\epsilon}(G/N)\rtimes_{\hat\epsilon|}(H/N)
\ar[u]^\cong
&
A\rtimes_\alpha G\rtimes_{\epsilon|}(G/H)
\ar[u]^\cong
\\
((A\rtimes_{\alpha}G
   \rtimes_{\epsilon}(G/N))\otimes C_0(G/H))
   \underset{\hat\epsilon\otimes\lt}{\rtimes}(G/N)
\ar[u]^{X_{H/N}^{G/N}(\hat\epsilon)}
&
A\rtimes_\alpha G\rtimes_{\epsilon}(G/N)
 \rtimes_{\hat\epsilon}(G/N)\rtimes_{\hat{\hat\epsilon}|}(G/H).
\ar[u]_{\hskip 0cm K(\epsilon)\rtimes(G/H)}
\ar[l]_-\cong
}
\end{equation}
of imprimitivity bimodules commutes, where all the isomorphisms are
canonical.
\end{thm}

We now show that Theorem~\ref{thm-YZ} follows from
Theorem~\ref{thm-KXZ}. 
See Section~\ref{next} for the proof of
Theorem~\ref{thm-KXZ}.  

\begin{proof}[Proof of Theorem~\ref{thm-YZ}]
Insert the arrow 
\begin{equation*}
\xymatrix@C=110pt{
A\rtimes_\alpha G\rtimes_{\epsilon}(G/N)\rtimes_{\hat\epsilon|}(H/N)
\ar[r]^-{Y_{G/H}^{G/N}(\epsilon)}
&
A\rtimes_\alpha G\rtimes_{\epsilon|}(G/H)
}
\end{equation*}
into the middle of commutative
diagram~\eqref{eq-KXZ} to create an upper square and a
lower square.  
Applying Corollary~6.4 of \cite{KQ-MI} to the maximal coaction
$\epsilon=\hat\alpha|$ of $G/N$ on $A\rtimes_{\alpha}G$ 
shows that the lower square commutes; since all arrows are invertible,
it follows that the upper square --- which is precisely~\eqref{YZ} ---
commutes as well.
\end{proof}

\medskip

In the remainder of this section, we prepare for the proof of
Theorem~\ref{thm-KXZ} by identifying 
%$X_{H/N}^{G/N}(\hat\epsilon)$ 
%and 
%$K(\epsilon)\rtimes_{\epsilon_K|}(G/H)$ 
the three bimodules in~\eqref{eq-KXZ} 
with symmetric-imprimitivity bimodules.
We retain the notation and hypotheses used thus far in this section
(but we will carefully note situations where in fact $H$ need not 
be normal in $G$).

\subsection{Realising $X_{H/N}^{G/N}(\hat\epsilon)$ as a 
symmetric-imprimitivity bimodule}
It is well-known how to use 
the symmetric imprimitivity theorem to derive 
Green's imprimitivity theorem \cite[Proposition~3]{green}.
It turns out that, for 
the action $\beta$ of $G/N$ on 
$(A\otimes C_0(G/N))\rtimes_{\alpha\otimes\lt} G)$ 
as at~\eqref{new-beta} and the subgroup $H/N\subseteq G/N$,
the symmetric imprimitivity theorem can produce the Green bimodule 
from a somewhat different set-up.
In this subsection, $H$ need not be normal in $G$.

First note that the identity map on $C_c(H/N\times G\times G/N, A)$ extends
to an isomorphism
\[
i:\big((A\otimes C_0(G/N)) \rtimes_{\alpha\otimes\lt}G\big)\rtimes_{\beta|}(H/N)\to (A\otimes C_0(G/N))\rtimes_\gamma (H/N\times G),\]
where $\gamma = (\id\times\alpha)\otimes(\rt\times\lt)
= (\id\otimes\rt)\times(\alpha\otimes\lt)$. 
The map 
\[\iota: C_c(G/N\times G/H\times G\times G/N,A)\to C_c(G/N\times G\times G/N\times G/H, A)\] 
defined by $\iota(g)(tN,s,rN,uH)=g(tN,uH,s,rN)$ extends to an isomorphism
\[\iota:((A\otimes C_0(G/N)\rtimes_{\alpha\otimes\lt} G)
\otimes C_0(G/H))\rtimes_{\beta\otimes\lt}(G/N)\to  (A\otimes C_0(G/N)\otimes C_0(G/H))
\rtimes_\eta(G/N\times G),\]
where $\eta$ is the action
$(\id\times\alpha)\otimes(\rt\times\lt)\otimes(\lt\times\id)
= (\id\otimes\rt\otimes\lt)\times(\alpha\otimes\lt\otimes\id)$.  

Now consider the symmetric imprimitivity
data $({}_KP_L,A,\sigma,\eta)$ defined as follows:
\begin{equation}\label{P-def}
\xymatrix{
& P = G/N\times G\times G/N & \\
K=G/N\times G
\ar[ru]^{(tN,s)\cdot(rN,u,vN) = (trN,su,tvN)\hspace{1.5cm}}
\ar[dr]_{\sigma=\id\times\alpha}&&
L=H/N\times G
\ar[lu]_{\hspace{1.5cm}(rN,u,vN)\cdot(hN,y) = (ryN,uy,vhN)}
\ar[dl]^{\eta=\id}\\
& A & }
\end{equation}
The symmetric imprimitivity theorem gives an 
$(\Ind_L^P\eta \rtimes_{\sigma\otimes\lt}K)-
(\Ind_K^P\sigma \rtimes_{\eta\otimes\rt}L)$
imprimitivity bimodule $W(P)=W({}_KP_{L}, A,\sigma,\eta)$.  

\begin{prop}\label{prop-green}
Suppose $\alpha\colon G\to\Aut A$ is a continuous action of a locally compact
group $G$ by automorphisms of a $C^*$-algebra $A$, and suppose $N$ and $H$
are closed subgroups of $G$ with $N$ normal in $G$ and  $N\subseteq H$.
Let $X_{H/N}^{G/N}(\beta)$ be the Green
bimodule associated to the action $\beta$ of $G/N$ on 
$(A\otimes C_0(G/N))\rtimes_{\alpha\otimes\lt}G$
described at~\eqref{new-beta},
%\[
%\big(\big((A\otimes C_0(G/N)\rtimes_{\alpha\otimes\lt} G)
%\otimes C_0(G/H)\big)\rtimes_{\beta\otimes\lt}(G/N)\big)
%-\big(\big((A\otimes C_0(G/N)) \rtimes_{\alpha\otimes\lt}G\big)
%\rtimes_{\beta|}(H/N)\big)
%\] 
and let $W(P)$ be the imprimitivity bimodule associated to the
symmetric imprimitivity data $({}_KP_L, A,\sigma,\eta)$
described at~\eqref{P-def}. 
Then there are \textup(non-canonical\textup)
equivariant isomorphisms 
$\Gamma\colon A\otimes C_0(G/N)\otimes C_0(G/H)\to \Ind_L^P\eta$ 
and 
$\Upsilon\colon A\otimes C_0(G/N)\to \Ind_K^P\sigma$
such that the diagram
\begin{equation}
\label{X}
\xymatrix{
{(((A\otimes C_0(G/N))\underset{\alpha\otimes\lt}{\rtimes}G)
     \otimes
C_0(G/H))\underset{\beta\otimes\lt}{\rtimes}(G/N)\hspace{.1cm}}
\ar[r]^-{X_{H/N}^{G/N}(\beta)}
\ar[d]_{(\Gamma\rtimes K)\circ\iota}^{\cong}
&
{\hspace{.1cm}((A\otimes C_0(G/N))\underset{\alpha\otimes\lt}{\rtimes}G)
    \underset{\beta|}{\rtimes}(H/N)}
\ar[d]^-{(\Upsilon\rtimes L)\circ i}_-{\cong}
\\
{\Ind_L^P\eta\rtimes_{\sigma\otimes\lt}K}
\ar[r]^-{W(P)}
&
{\Ind_K^P\sigma\rtimes_{\eta\otimes\rt}L}
}
\end{equation}
%\begin{equation}
%\label{Xold}
%\xymatrix{
%{(((A\otimes C_0(G/N))\rtimes_{\alpha\otimes\lt}G)\otimes C_0(G/H))\rtimes_{\beta\otimes\lt}(G/N)}
%\ar[d]_{X_{H/N}^{G/N}}
%\ar[rr]^{\hskip3cm(\Gamma\rtimes K)\circ\iota}_{\hskip3cm\cong}
%&&
%{\Ind_L^P\eta\rtimes_{\sigma\otimes\lt}K}
%\ar[d]^{W({}_KP_L)}
%\\
%{((A\otimes C_0(G/N))\rtimes_{\alpha\otimes\lt}G)\rtimes_{\beta|}(G/N)}
%\ar[rr]_-{(\Upsilon\rtimes L)\circ i}^-{\cong}
%&&
%{\Ind_K^P\sigma\rtimes_{\eta\otimes\rt}L}
%}
%\end{equation}
of imprimitivity bimodules commutes.
\end{prop}

\begin{proof}
For $f\in A\otimes C_0(G/N)\cong C_0(G/N, A)$ and $(rN, u, vN)\in  G/N\times G\times G/N$ define
\begin{equation}\label{eq-Upsilon}
\Upsilon(f)(rN,u,vN) = \alpha_u(f(r^{-1}vN)).
\end{equation}
Then, for $(tN,s)\in K=G/N\times G$,
\begin{align*}
\Upsilon(f)((tN,s)\cdot(rN,u,vN))&=\Upsilon(f)(trN,su, tvN)=\alpha_{su}(f(r^{-1}vN))\\
&=\alpha_s(\Upsilon(f)(rN,u, vN))=\sigma_{(tN,s)}(\Upsilon(f)(rN,u, vN)),
\end{align*}
so  $\Upsilon$ maps $A\otimes C_0(G/N)$ into $\Ind_K^P\sigma$.  
It is straightforward to check that $\Upsilon$ is invertible, 
with inverse given by
$\Upsilon^{-1}(g)(tN) = g(N,e,tN)$.
For $(hN,y)\in L=N/H\times G$,
\begin{align*}
\Upsilon(\gamma_{(hN,y)}(f))(rN, u, vN)
&=\alpha_u(\gamma_{(hN,y)}(f)(r^{-1}vN))\\
&=\alpha_{uy}(f(y^{-1}r^{-1}vhN))\\
&=\Upsilon(f)(ryN, uy, vhN)\\
&=\Upsilon(f)((rN, u, vN)\cdot(hN, y))\\
&=(\eta\otimes\rt)_{(hN,y)}(\Upsilon(f))(rN, u, vN),
\end{align*}
so $\Upsilon$ is a  $\gamma$ -- $(\eta\otimes\rt)$ equivariant isomorphism and
induces an isomorphism
\[\Upsilon\rtimes L: (A\otimes C_0(G/N))\rtimes_\gamma L\to \Ind_K^P\sigma\rtimes_{\eta\otimes\rt}L\]
of the crossed products.

Similarly, the map $\Gamma\colon A\otimes C_0(G/N)\otimes C_0(G/H)
\to \Ind_L^P\eta$
defined by
\begin{equation}\label{eq-Gamma}
\Gamma(f)(rN,u,vN) = f(ur^{-1}N,vH)
\end{equation}
is an $\epsilon$ -- $(\sigma\otimes\lt)$ equivariant
isomorphism with inverse given by
$\Gamma^{-1}(g)(vN,rH) = g(v^{-1}N,e,rN)$. So $\Gamma$ induces an isomorphism
\[
\Gamma\rtimes K:(A\otimes C_0(G/N)\otimes C_0(G/H))
\rtimes_\epsilon K\to\Ind_L^P\eta\rtimes_{\sigma\otimes\lt} K.
\]

Let $\Psi:C_c(G/N\times G\times G/N, A)\to C_c(G/N\times G\times G/N, A)$
be the map \[\Psi(rN, u, vN)=f(vN, u, ur^{-1}N)\Delta_G(u)^{1/2}.\]
We will show that the triple $((\Gamma\rtimes K)\circ\iota,
\Psi,(\Upsilon\rtimes L)\circ i)$ extends to an imprimitivity bimodule
isomorphism of $X_{H/N}^{G/N}(\beta)$ onto $W(P)$.  We may view both
$X_{H/N}^{G/N}(\beta)$ and  $W(P)$ as  completions of $C_c(G/N\times
G\times G/N, A)$, so  $\Psi$  clearly has dense range. It therefore
suffices to show, for  $x,y\in C_c(G/N\times G\times G/N,A)\subseteq
X_{H/N}^{G/N}(\beta)$ and  $f\in C_c(H/N\times G\times G/N,A)$, that
\begin{enumerate}
\item
$\Psi(x\cdot f)=\Psi(x)\cdot ((\Upsilon\rtimes L)\circ i(f))$; and
\item
$(\Gamma\rtimes K)\circ \iota({}_*\langle x\,,\, y\rangle)={}_*\langle\Psi(x)\,,\, \Psi(y)\rangle$.
\end{enumerate}
(For then (2) implies \[\|\Psi(g\cdot x)-((\Gamma\rtimes K)\circ\iota(g))\cdot\Psi(x)\|^2=0\] for $g\in C_c(G/N\times G/H\times G\times G/N, A)$, 
and this together with (1), (2), and denseness gives the other inner product condition.)

So let $x,y$ and $f$ be as above. Using the formula for the right action in Green's bimodule from \cite[Equation B.5]{BE} we have:
\begin{align*}
\Psi&(x\cdot f)(rN, u, vN)
=\Delta_G(u)^{1/2}(x\cdot f)(vN, u, ur^{-1}N)\\
&=\Delta_G(u)^{1/2}\int_{H/N}x(vhN,\cdot,\cdot)\beta_{vhN}(f(h^{-1}N,\cdot,\cdot))\Delta_{H/N}(hN)^{-1/2}\, d(hN)\, (u,ur^{-1}N)\\
&=\Delta_G(u)^{1/2}\int_{H/N}\int_G x(vhN, s,\cdot)(\alpha\otimes\lt)_s(\beta_{vhN}(f(h^{-1}N, s^{-1}u, \cdot)))\\
&\hskip 8cm \Delta_{H/N}(hN)^{-1/2}\, d(hN)\, ds\, (ur^{-1}N)\\
&=\Delta_G(u)^{1/2}\int_{H/N}\int_G x(vhN, s,ur^{-1}N)\alpha_s(f(h^{-1}N, s^{-1}u, s^{-1}ur^{-1}vhN))\\
&\hskip 9cm \Delta_{H/N}(hN)^{-1/2}\, d(hN)\, ds.
\end{align*}
Using the formula for the right action on $W(P)$ from \cite[Equation~B.2]{BE} we have
\begin{align*}
\big(\Psi&(x)\cdot((\Upsilon\rtimes L)\circ i(f))\big) (rN, u, vN)\\
&=\int_{H/N\times G}\eta_{(hN,t)}(\Psi(x)((rN, u, vN)\cdot(hN,t))\\
&\hskip 2cm\Upsilon\rtimes L(f)((hN,t)^{-1},(rN,u,vN)\cdot(hN,t)))\Delta_{H/N\times G}((hN,t))^{-1/2}\, d(hN,t)\\
&=\int_{H/N}\int_G\Psi(x)(rtN, ut, vhN)\Upsilon(f(h^{-1}N, t^{-1}))(rtN, utN, vhN)\\
&\hskip 8cm \Delta_{H/N}(hN)\Delta_G(t)^{-1/2}\, d(hN)\, dt\\
&=\Delta_G(u)^{1/2}\int_{H/N}\int_G x(vhN, ut, ur^{-1}N)\alpha_{ut}(f(h^{-1}N, t^{-1},t^{-1}r^{-1}vhN))
\\
&\hskip 10cm\Delta_{H/N}(hN)^{-1/2}\, d(hN)\, dt,
\end{align*}
which equals $\Psi(x\cdot f)(rN, u, vN)$ by the change of variable $s=ut$.
Also,
\begin{align*}
(\Gamma&\rtimes K)\circ\iota({}_*\langle x\,,\, y\rangle)((tN,s),(rN,u,vN))
=\Gamma(\iota({}_*\langle x\,,\, y\rangle)(tN,s,\cdot,\cdot))(rN,u,vN)\\
&={}_*\langle x\,,\, y\rangle(tN,vH, s,ur^{-1}N)\\
&=\Delta_{G/N}(tN)^{-1/2}\int_{H/N}x(vhN,\cdot,\cdot)\beta_{tN}(y(t^{-1}vhN,\cdot,\cdot)^*)\, d(hN) (s, ur^{-1}N)\\
\intertext{(note that the product
$x(vhN,\cdot,\cdot)\beta_{tN}(y(t^{-1}vhN,\cdot,\cdot)^*)$ is
convolution in $(A\otimes C_0(G/N))\rtimes_{\alpha\otimes\lt}G$)}
&=\Delta_{G/N}(tN)^{-1/2}\int_{H/N}\int_G x(vhN,w,\cdot)(\alpha\otimes\lt)_w (\beta_{tN}(y(t^{-1}vhN,\cdot,\cdot)^*))(w^{-1}s)\\
&\hskip 11cm \, dw\, d(hN) (ur^{-1}N)\\
&=\Delta_{G/N}(tN)^{-1/2}\int_{H/N}\int_G x(vhN,w,ur^{-1}N)(\alpha\otimes\lt)_w (y(t^{-1}vhN,\cdot,\cdot)^*)(w^{-1}s,ur^{-1}tN)\\
&\hskip 13cm \, dw\, d(hN) \\
\intertext{(note that the involution $y(t^{-1}vhN,\cdot,\cdot)^*$ is in
$(A\otimes C_0(G/N))\rtimes_{\alpha\otimes\lt}G$)}
&=\Delta_{G/N}(tN)^{-1/2}\int_{H/N}\int_G x(vhN,w,\cdot) (\alpha\otimes\lt)_{ww^{-1}s}(y(t^{-1}vhN,s^{-1}w,\cdot)^*)(ur^{-1}tN)\\
&\hskip 12cm\Delta_G(s^{-1}w)\, dw\, d(hN)\\
&=\Delta_{G/N}(tN)^{-1/2}\int_{H/N}\int_G x(vhN,w,ur^{-1}N)
\alpha_{s}(y(t^{-1}vhN,s^{-1}w,s^{-1}ur^{-1}tN)^*)\\*
&\hskip 12cm\Delta_G(s^{-1}w)\, dw\, d(hN),\\
&=\Delta_{G/N}(tN)^{-1/2}\int_{H/N}\int_G x(vhN, uw,
ur^{-1}N)\alpha_s(y(t^{-1}vh, s^{-1}uw, s^{-1}ur^{-1}tN)^*)) \\*
&\hskip11cm
\Delta_G(s^{-1}uw)\, d(hN) \, dw\\
&=\int_{H/N}\int_G \Psi(x)(rwN, uw, vhN)\alpha_s(\Psi(y)(t^{-1}rwN,
s^{-1}uw, t^{-1}vhN)^*)\\*
&\hskip8cm \Delta_{G/N}(tN)^{-1/2}\Delta_G(s)^{-1/2}\, d(hN) \, dw\\
&=\int_{H/N\times G}\eta_{(hN, w)}\big( \Psi(x)((rN,u,vN)\cdot (hN,w))\\*
&\hskip1cm \sigma_{(tN,s)}(\Psi(y)((tN,s)^{-1}(rN,u,vN)\cdot(hN,w))^*) \big)\Delta_{G/N\times G}((tN, s))^{-1/2}\, d(hN)\, dw\\
&={}_*\langle \Psi(x)\,,\, \Psi(y)\rangle((tN,s),(rN,u,vN)).
\end{align*}
\end{proof}

%----------------------------------------------------------------------
\subsection{The definition of $Z_{G/H}^{G/N}(\alpha)$}
\label{Q-sec}

In \cite{aHKRW-EP},  $Z_{G/H}^{G/N}(\alpha)$ is defined using
symmetric imprimitivity data $({}_LQ_G,A,\zeta,\tau)$ as follows:
\begin{equation}\label{Q-def}
\xymatrix{
& Q = G/N\times G & \\
L=H/N\times G
\ar[ru]^{(hN,x)\cdot(wN,z) = (hwN,xz)\hspace{1cm}}
\ar[dr]_{\zeta=\id\times\alpha}&&
\hspace{1cm}G\hspace{1cm}
\ar[lu]_{\hspace{1cm}(wN,z)\cdot y = (wyN,zy)}
\ar[dl]^{\tau=\id}\\
& A & }
\end{equation}
(Here again,  $H$ need not be normal in $G$.)
More precisely, 
the symmetric imprimitivity theorem gives an $(\Ind_G^Q\tau
\rtimes_{\zeta\otimes\lt}L)$--$(\Ind_L^Q\zeta \rtimes_{\tau\otimes\rt}G)$
imprimitivity bimodule $W(Q)=W({}_LQ_G, A,\zeta,\tau)$.
The map $\Omega\colon A\otimes C_0(G/H)\to \Ind_L^Q\zeta$
defined by
\begin{equation}\label{eq-Omega}
\Omega(f)(wN,z) = \alpha_z(f(w^{-1}H))
\end{equation}
is an $(\alpha\otimes\lt)$--$(\tau\otimes\rt)$ equivariant isomorphism
with inverse given by $\Omega^{-1}(g)(tH) = g(t^{-1}N,e)$ and hence induces an isomorphism $\Omega\rtimes G$ of $(A\otimes C_0(G/H))\rtimes_{\alpha\otimes\lt}G$ onto $\Ind_L^Q\zeta\rtimes_{\tau\otimes\rt} G$.
The map $\Theta\colon A\otimes C_0(G/N)\to \Ind_G^Q\tau$
defined by
\begin{equation}\label{eq-Theta}
\Theta(f)(wN,z) = f(zw^{-1}N)\quad\quad((wN,z)\in Q=G/N\times G)
\end{equation}
is a $\gamma$ -- $(\zeta\otimes\lt)$
equivariant isomorphism with inverse given by 
$\Theta^{-1}(g)(rN) = g(r^{-1}N,e) = g(N,r)$.
So $\Theta$  induces an isomorphism $\Theta\rtimes L$  of $(A\otimes
C_0(G/N))\rtimes_\gamma L$ onto
$\Ind_G^Q\tau\rtimes_{\zeta\otimes\lt}L$.  
The imprimitivity bimodule
$Z_{G/H}^{G/N}(\alpha)$ is then defined by 
requiring the diagram
\begin{equation}
\label{eq-W(Q)}
\xymatrix{
{\big((A\otimes C_0(G/N))\rtimes_{\alpha\otimes\lt}G\big)\rtimes_{\beta|}(H/N)}
\ar[rr]^-{Z_{G/H}^{G/N}(\alpha)}
\ar[d]_{(\Theta\rtimes L)\circ i}^{\cong}
&&
{(A\otimes C_0(G/H))\rtimes_{\alpha\otimes\lt}G.}
\ar[d]^{\Omega\rtimes G}_{\cong}
\\
{\Ind_G^Q\tau\rtimes_{\zeta\otimes\lt}L}
\ar[rr]^-{W(Q)}
&&
{\Ind_L^Q\zeta\rtimes_{\tau\otimes\rt}G}
}
\end{equation}
to commute.

%----------------------------------------------------------------------

\subsection{Realising $K(\epsilon)\rtimes_{\epsilon_K|}(G/H)$ as a
symmetric-imprimitivity bimodule}
The most difficult bimodule in Theorem~\ref{thm-KXZ}
is the crossed-product Katayama bimodule
$K(\epsilon)\rtimes_{\epsilon_K|}(G/H)$.
The difficulty arises partly because of the coaction crossed-product,
and partly because the Katayama bimodule is inherently spatial.
We were able to obtain this realisation
by looking at a set-up which should implement $K(\epsilon)$,
and then adding $G/H$ with the appropriate group actions.

Consider the symmetric imprimitivity data
\begin{equation}\label{R-def}
\xymatrix{
 & R=G/N\times G\times G/H&  \\
K=G/N\times G\ar[ru]^{(tN,s)\cdot(rN,u,vH) = (trN,su,tvH)\hspace{1.2cm}}
 \ar[rd]_{\sigma=\id\times\alpha}&  & 
\hspace{1cm}G\hspace{1cm}\ar[lu]_{\hspace{1cm}(rN,u,vH)\cdot y = (ryN,uy,vH)}
 \ar[ld]^{\tau=\id}\\
 & A&  }
\end{equation}
(Note that $K$ and $\sigma$ are the same as for $P$, 
and $\tau$ is the same as for $Q$.)
The symmetric imprimitivity theorem gives an
$(\Ind_G^R\tau \rtimes_{\sigma\otimes\lt}K)$--$(\Ind_K^R\sigma
\rtimes_{\tau\otimes\rt}G)$ imprimitivity bimodule 
$W(R)= W({}_KR_G,A,\sigma,\tau)$.

\begin{prop}\label{new-R}
Suppose $\alpha :G\to\Aut A$ is a continuous action of a locally compact
group $G$ by automorphisms of a $C^*$-algebra $A$, and suppose $N$ and $H$
are closed normal subgroups of $G$ with $N\subseteq H$. 
Let $K(\epsilon)$ be the Katayama bimodule 
as defined at~\eqref{eq-maximal-iso} associated
to the maximal coaction $\epsilon=\hat\alpha|_{G/N}$
of $G/N$ on $A\rtimes_\alpha G$,
and let $W(R)$ be the bimodule associated to the symmetric imprimitivity
data $({}_KR_G,A,\sigma,\tau)$ described at~\eqref{R-def}.
Then there exist \textup(non-canonical\textup)
equivariant isomorphisms
$\Lambda\colon A\otimes C_0(G/N)\otimes C_0(G/H)\to\Ind_G^R\tau$
and
$\Xi\colon A\otimes C_0(G/H)\to\Ind_K^R\sigma$
such that the diagram
\begin{equation}\label{KR}
\xymatrix{
A\rtimes_\alpha G \rtimes_{\epsilon}(G/N)\rtimes_{\hat\epsilon}(G/N)
\rtimes_{\hat{\hat\epsilon}|}(G/H)
\ar[d]_\cong
\ar[r]^-{K(\epsilon)\rtimes(G/H)}
&
A\rtimes_{\alpha}G \rtimes_{\epsilon|}(G/H)
\ar[d]^\cong
\\
(((A\otimes C_0(G/N))\underset{\alpha\otimes\lt}{\rtimes}G)
   \otimes C_0(G/H))\underset{\beta\otimes\lt}{\rtimes}(G/N)
\ar[d]_{(\Lambda\rtimes K)\circ\iota}
&
(A\otimes C_0(G/H))\underset{\alpha\otimes\lt}{\rtimes}G
\ar[d]^{\Xi\rtimes G}
\\
\Ind_G^R\tau\rtimes_{\sigma\otimes\lt}K
\ar[r]^-{W(R)}
&
\Ind_K^R\sigma\rtimes_{\tau\otimes\rt}G
}
\end{equation}
of imprimitivity bimodules commutes, 
where the unnamed isomorphisms are the canonical ones.
\end{prop}

\begin{proof}
The map
$\Lambda\colon A\otimes C_0(G/N)\otimes C_0(G/H)\to\Ind_G^R\tau$
defined by
\begin{equation}\label{eq-Lambda}
\Lambda(f)(rN,u,vH) = f(ur^{-1}N,vH)
\end{equation}
is an $\epsilon - (\sigma\otimes\lt)$ equivariant
isomorphism with inverse given by
$\Lambda^{-1}(g)(rN,vH) = g(r^{-1}N,e,vH)$
%(This can be checked directly, but it will follow
%in Section~\ref{isoms}
%from the identity
%$\xLambda = \Gamma^{-1}\circ\phi_\tau\circ\psi_\tau$.)
So $\Lambda$ induces an isomorphism 
$\Lambda\rtimes K$ of 
$\big((A\otimes C_0(G/N)\otimes C_0(G/H)\big)\rtimes_\epsilon K$
onto 
$\Ind_G^R\tau\rtimes_{\sigma\otimes\lt}K$.
The map $\Xi\colon A\otimes C_0(G/H)\to\Ind_K^R\sigma$
defined by
\begin{equation}\label{eq-Xi}
\Xi(f)(rN,u,vH) = \alpha_u(f(r^{-1}vH))
\end{equation}
is an $(\alpha\otimes\lt)-(\tau\otimes\rt)$ equivariant isomorphism
with inverse given by 
$\Xi^{-1}(g)(wH) = g(N,e,wH)$.
%(This will follow from
%$\xXi = \Omega^{-1}\circ\phi_\sigma\circ\psi_\sigma$
%in Section~\ref{isoms}.)
So $\Xi$ also induces an isomorphism $\Xi\rtimes G$ of the crossed
products.  

We now define an imprimitivity bimodule $W$ to be $W(R)$ with the coefficient
algebras adjusted using these isomorphisms.  Thus, the following diagram
commutes by definition:
\begin{equation}
\label{eq-Wdefn}
\xymatrix{
{\big(((A\otimes C_0(G/N))\rtimes_{\alpha\otimes\lt}G)\otimes C_0(G/H)\big)\rtimes_{\beta\otimes\lt}(G/N)}
\ar[rr]^-{W}
\ar[d]_{(\Lambda\rtimes K)\circ\iota}^{\cong}
&&
{(A\otimes C_0(G/H))\rtimes_{\alpha\otimes\lt}G}
\ar[d]^{\Xi\rtimes G}_{\cong}
\\
{\Ind_G^R\tau\rtimes_{\sigma\otimes\lt}K}
\ar[rr]^-{W(R)}
&&
{\Ind_K^R\sigma\rtimes_{\tau\otimes\rt}G.}
}
\end{equation}
The formulas for the actions and inner products
of $W$ are as follows:
for $b\in C_c(G/N\times G/H\times G\times G/N,A)$,
$f,g\in C_c(G/N\times G\times G/H,A)\subseteq W$ 
and
$c\in C_c(G\times G/H,A)$,
\begin{align}\label{eq-Wequns}
\begin{split}
(b\cdot f)(rN,u,vH)
&= \int_{G/N}\int_G b(tN,vH,s,rN)
\alpha_s(f(s^{-1}rtN,s^{-1}u,t^{-1}vH))\\
&\phantom{xxxxxxxxxxxxxxxxxxxxxxxxxxx}
\Delta_G(s)^{1/2}\Delta_{G/N}(tN)^{1/2}\, \,ds \,d(tN)\\
(f\cdot c)(rN,u,vH)
&= \int_G f(rN,uy,vH) \alpha_{uy}(c(y^{-1},y^{-1}u^{-1}rvH))
\Delta_G(y)^{-1/2} \, dy\\
{}_*\langle f,g\rangle(tN,vH,s,rN)
&= \Delta_G(s)^{-1/2}\Delta_{G/N}(tN)^{-1/2}\int_G f(rN,y,vH)\\
&\hskip 5 cm\alpha_s(g(s^{-1}rtN,s^{-1}y,t^{-1}vH)^*)
\,dy
\\
\langle f,g\rangle_*(y,wH)
&= \Delta_G(y)^{-1/2}\int_G\int_{G/N} \alpha_s( f(s^{-1}tN,s^{-1},t^{-1}wH)^*\\
&\hskip5 cm g(s^{-1}tN,s^{-1}y,t^{-1}wH))
\,d(tN)\,ds.
\end{split}
\end{align}
\renewcommand{\qedsymbol}{}
\end{proof}

To complete the proof of Proposition~\ref{new-R}, 
it suffices to show that $W$ and
$K(\epsilon)\rtimes_{\epsilon_K|}(G/H)$ are isomorphic, 
modulo the canonical isomorphisms of the coefficient algebras.
This will involve a spatial argument.  
Recall
from \cite[Definition~2.1]{ER-MI} that a representation of an
$A$--$B$ imprimitivity bimodule $X$ on a pair of Hilbert spaces
$(\H_\llt,\H_\rrt)$ is a triple $(\mu_\llt,\mu,\mu_\rrt)$ consisting of
non-degenerate representations $\mu_\llt:A\to B(\H_\llt)$, $\mu_\rrt:B\to
B(\H_\rrt)$, and a linear map $\mu:X\to B(\H_\rrt,\H_\llt)$ such that,
for all $x,y\in X$, $a\in A$ and $b\in B$, \begin{enumerate}
\item $\mu(x)^*\mu(y)=\mu_\rrt(\langle x\,,\, y\rangle_B)$ and
$\mu(x)\mu(y)^*=\mu_\llt({}_A\langle x\,,\, y\rangle)$ and \item
$\mu(a\cdot x\cdot b)=\mu_\llt(a)\mu(x)\mu_\rrt(b)$.  \end{enumerate} The
representation $(\mu_\llt,\mu,\mu_\rrt)$ is faithful if either $\mu_\llt$
or $\mu_\rrt$ is isometric (for then $\mu$ is also isometric).

\begin{lem}\label{lem-coaction-x}
Let $(\mu_\llt,\mu,\mu_\rrt)$ be a faithful representation of
an imprimitivity bimodule ${}_AX_B$ on a pair of Hilbert spaces
$(\H_\llt,\H_\rrt)$. Let ${}_{\delta_A}\delta_{\delta_B}$ be a full coaction of
$G$ on ${}_AX_B$, so that $X\rtimes_\delta G$ is an
$(A\rtimes_{\delta_A}G)$--$(B\rtimes_{\delta_B}G)$ imprimitivity bimodule.   
Let $\mu_\llt\rtimes
G$ and $\mu_\rrt\rtimes G$ be the regular representations of
$A\rtimes_{\delta_A}G$ and $B\rtimes_{\delta_B}G$ induced from $\mu_\llt$
and $\mu_\rrt$, respectively, and let
\[
\mu\rtimes G:=(\mu\otimes\lambda)\circ\delta\rtimes(1\otimes
M):X\rtimes_\delta G\to B\big(\H_\rrt\otimes L^2(G),\H_\llt\otimes L^2(G)\big).
\]
Then $(\mu_\llt\rtimes G,\mu\rtimes G,\mu_\rrt\rtimes G)$ is a faithful
representation of $X\rtimes_\delta G$.
\end{lem}

\begin{proof}
This is essentially Theorem~3.2 of \cite{ER-MI}. However, since it is
proved there for reduced coactions, we outline an alternative proof based
on results in \cite{BE}. The representations $(\mu_\llt,\mu,\mu_\rrt)$
combine to give a faithful representation $L(\mu)$ of the linking
algebra $L(X)$ as bounded operators on $\H_{\llt}\oplus\H_{\rrt}$. As
in \cite[Chapter~3, \S1.2]{BE}, the coactions $\mu_\llt$, 
$\mu$, and $\mu_\rrt$ combine to give a
coaction $\nu$ of $G$ on $L(X)$, and $L(X)\rtimes_\nu G$ is canonically
isomorphic to $L(X\rtimes_\delta G)$ \cite[Lemma~3.10]{BE}. The regular
representation $L(\mu)\rtimes G$ of $L(X)\rtimes_\nu G$ on
\[
(\H_{\llt}\oplus\H_{\rrt})\otimes L^2(G)
=(\H_{\llt}\otimes L^2(G))\oplus(\H_{\rrt}\otimes L^2(G))
\]
is faithful by \cite[Corollary~A.59]{BE}. Since
\[
L(\mu)\rtimes G=((L(\mu)\otimes \lambda)\circ \nu)\rtimes (1\otimes M)
\]
restricts to the regular representations on the corners of
$L(X\rtimes_\delta G)$, we deduce that $\mu\rtimes G$ is faithful too.
\end{proof}

%\begin{thm}\label{thm-KxGmodM} Let $\alpha :G\to\Aut A$ be a continuous
%action of a locally compact group $G$ by automorphisms of a
%$C^*$-algebra $A$. Also let $N$ and $H$ be closed normal subgroups of $G$ with $N\subseteq H$.
%Let $Z$ be the 
%\[
%\big(\big(((A\otimes C_0(G/N))\rtimes_{\alpha\otimes\lt}G)\otimes C_0(G/H)\big)\rtimes_{\beta\otimes\lt}(G/N)\big)-\big((A\otimes C_0(G/H))\rtimes_{\alpha\otimes\lt}G\big)
%\] 
%imprimitivity bimodule $Z$  defined at \eqref{eq-Zdefn}. Set $\delta=\hat\alpha|_{G/N}$. Then $Z$ and the \[\big(A\rtimes_\alpha G\rtimes_{\delta}(G/N)\rtimes_{\hat\delta}(G/N)\rtimes_{\hat{\hat\delta}|}(G/H)\big)-\big(A\rtimes_\alpha G\rtimes_{\hat\alpha|} (G/H)\big)\] imprimitivity bimodule
%$K(\delta)\rtimes_{\delta_{K|}}(G/H)$ are isomorphic.
%\end{thm}

\begin{proof}[Conclusion of the Proof of Proposition~\ref{new-R}]

Let $(\pi,U)$ be a faithful covariant representation of $(A,G,\alpha)$
on a Hilbert space $\H$.  The idea of the proof is to find faithful
representations $(\nu_\llt, \nu,\nu_\rrt)$  and $\big(\mu_\llt\rtimes
(G/H),\mu\rtimes (G/H),\mu_\rrt\rtimes (G/H)\big)$ of $W$ and
$K(\epsilon)\rtimes_{\epsilon_{K|}}(G/H)$  on
\[
\big(\H\otimes L^2(G/N)\otimes L^2(G/H), \H\otimes L^2(G/H)\big)
\]
such that the ranges of $\nu_\rrt$ and $\mu_\rrt\rtimes (G/H)$ coincide. We
will then argue that a dense subset of the range of $\nu$ is contained
in the range of $\mu\rtimes G/H$.  Thus $W$ is isomorphic to a
closed submodule of $K(\epsilon)\rtimes_{\epsilon_{K|}}(G/H)$
on which the right inner product is full, and  it then follows
from the Rieffel correspondence that  $W$ and 
$K(\epsilon)\rtimes_{\epsilon_{K|}}(G/H)$ are isomorphic.

The  representation
\[\nu_\rrt:=(\pi\otimes M^{G/H})\rtimes(U\otimes\lambda^{G/H}):(A\otimes
C_0(G/H))\rtimes_{\alpha\otimes\lt}G\to B(\H\otimes L^2(G/H))\] is
faithful; for future use, note that it is given on the pieces
$A$, $C^*(G)$ and $C_0(G/H)$ by
$\pi\otimes 1$, $U\otimes\lambda^{G/H}$ and $1\otimes M^{G/H}$, respectively.  Let
\[
\nu_\llt:\big(((A\otimes C_0(G/N))\rtimes_{\alpha\otimes\lt}G)\otimes C_0(G/H)\big)\rtimes_{\beta\otimes\lt}G/N\to B\big(\H\otimes L^2(G/N)\otimes L^2(G/H)\big)
\]
be the representation \[
\nu_\llt:=\big((\pi\otimes M^{G/N}\rtimes U\otimes\lambda^{G/N})\otimes
M^{G/N}\big)\rtimes 1\otimes \rho\otimes\lambda^{G/H};
\]
it is given on the pieces
$A$, $C^*(G)$, $C_0(G/N)$, $C^*(G/N)$ and $C_0(G/H))$ by
$\pi\otimes1\otimes1$, $U\otimes\lambda^{G/N}\otimes 1$, $1\otimes
M^{G/N}\otimes 1$, $1\otimes\rho\otimes\lambda^{G/H}$ and $1\otimes1\otimes M^{G/H}$, respectively.
Next, we claim that, for
fixed $z\in C_c(G/N\times G\times G/H,A)\subseteq W$ and every $\xi\in L^2(G/H,\H)$,
 the map $\nu(z)\xi$ of
$G/N\times G/H$ into $\H$ given by
\begin{equation}
\label{eq-nu}
(\nu(z)\xi)(rN,vH)
= \int_G \pi(z(rN,y,vH))U_y\xi(y^{-1}rvH)
  \Delta_G(y)^{-1/2}\, dy
\end{equation}
is an element of $L^2(G/N\times G/H,\H)\cong\H\otimes L^2(G/N)\otimes L^2(G/H)$.
For $w\in C_c(G/N\times G\times G/H,A)$
and $\eta\in L^2(G/H,\H)$, we have
\begin{align*}
&\int_{G/N\times G/H}\bigl( \nu(w)\eta(rN,vH)
\bigm| \nu(z)\xi(rN,vH) \bigr) \,d(rN,vH)\\
&= \int_{G/N}\int_{G/H}\int_G\int_G
\bigl( \pi(w(rN,x,vH)) U_x \eta(x^{-1}rvH)
\bigm| \pi(z(rN,y,vH))U_y \xi(y^{-1}rvH) \bigr) \\
&\phantom{\int_{G/H}\int_G\int_{G\times G/N}}
\Delta_G(xy)^{-1/2} \,dx\,dy\,d(vH)\,d(rN)\\
&\overset{(\dagger)}{=} \int_{G/N}\int_{G/H}\int_G\int_G\
\bigl( \pi(w(y^{-1}rN,y^{-1}x,r^{-1}vH)) U_{y^{-1}x}
\eta(x^{-1}vH)\\
&\phantom{\int_{G/H}\int_G\int_G\int_{G/N}}
\bigm| \pi(z(y^{-1}rN,y^{-1},r^{-1}vH))U_{y^{-1}}
\xi(vH) \bigr) \Delta_G(x)^{-1/2} \,dx\,dy\,d(vH)\,d(rN)\\
&= \int_{G/H}\int_G\int_{G}\int_{G/N}
\bigl( U_y\pi(z(y^{-1}rN,y^{-1},r^{-1}vH)^*
w(y^{-1}rN,y^{-1}x,r^{-1}vH)) U_{y^{-1}x} \eta(x^{-1}vH)\\
&\phantom{\int_{G/H}\int_G\int_{G\times G/N}}
\bigm| \xi(vH) \bigr) \Delta_G(x)^{-1/2} \,d(y,rN)\,dx\,d(vH)\\
&=  \int_{G/H}\int_G
\Big( \pi\big(\int_G\int_{G/N}\alpha_y(z(y^{-1}rN,y^{-1},r^{-1}vH)^*
w(y^{-1}rN,y^{-1}x,r^{-1}vH))\,d(y)\, d(rN)\big)\\
&\phantom{\int_{G/H}\int_G\int_{G\times G/N}}
U_x \eta(x^{-1}vH)\bigm| \xi(vH) \Big) \Delta_G(x)^{-1/2} \,dx\,d(vH)\\
&= \int_{G/H}\int_G \bigl( \pi(\langle z,w\rangle_*(x,vH)) U_x
\eta(x^{-1}vH) \bigm| \xi(vH) \bigr) \,dx\,d(vH)\\
&= \int_{G/H} \bigl( ((\pi\otimes M^{G/H})\rtimes(U\otimes\lambda^{G/H})
(\langle z,w\rangle_*)\eta)(vH) \bigm| \xi(vH) \bigr) \,d(vH)\\
&=\bigl( \nu_\rrt(\langle z,w\rangle_*)\eta
\bigm| \xi \bigr),
\end{align*}
where $(\cdot\mid\cdot)$ denotes the appropriate
Hilbert space inner product.  The change of variables at $(\dagger)$ is given by
$(vH,rN,x,y)\mapsto(r^{-1}vH,y^{-1}rN,y^{-1}x,y^{-1})$. In particular, this shows that
$\|\nu(z)\xi\|_2^2
=\bigl( \nu(z)\xi \bigm| \nu(z)\xi \bigr)
\leq \|z\|^2\|\xi\|^2$,
so $\nu(z)\xi\in L^2(G/N\times G/H,\H)$, and that the linear map $\xi\mapsto \nu(z)\xi$ is bounded.  Thus $\nu$, as defined at \eqref{eq-nu}, extends to a linear map
\[\nu:W\to B\big(\H\otimes L^2(G/H),\H\otimes L^2(G/N)\otimes L^2(G/H)\big).\]

We claim that $(\nu_\llt,\nu,\nu_\rrt)$ is a representation of $W$. We will prove
that, for $z,w\in C_c(G/N\times G\times G/H, A)\subseteq W$ and $b\in C_c( G/N\times G/H\times G\times G/N,A)$,
\begin{enumerate}
\item $\nu(z)^*\nu(w)=\nu_\rrt(\langle z\,,\, w\rangle_*)$ in $B(\H\otimes L^2(G/H))$;
\item $\nu(b\cdot z)=\nu_\llt(b)\nu(z)$; and
\item $\nu$ is non-degenerate in the sense that $\{\nu(z)\xi:z\in W,\xi\in \H\otimes L^2(G/H)\}$ is dense in $\H\otimes L^2(G/N)\otimes L^2(G/H)$.
\end{enumerate}
Then (1) implies that $\nu(z\cdot c)=\nu(z)\nu_\rrt(c)$ for all $c\in (A\otimes C_0(G/H))\rtimes_{\alpha\otimes\lt}G$, and the other inner product condition follows from this and (1)--(3).

To see that (1) holds, it suffices to see that
\[
\bigl( \nu(w)\eta \bigm| \nu(z)\xi \bigr)
= \bigl( \nu_\rrt(\langle z,w\rangle_*)\eta
\bigm| \xi \bigr)
\]
for all $z,w\in C_c(G/N\times G\times G/H,A)$
and $\xi,\eta\in L^2(G/H,\H)$, and this was done
in the calculation above which showed $\nu$ is well-defined.

It will be easiest to check (2) on the separate pieces of the algebra.
The piece
$A\otimes C_0(G/N)\otimes C_0(G/H)$ is represented
by $\pi\otimes M^{G/N}\otimes M^{G/H}$ and we deduce from \eqref{eq-Wequns} that $b\in C_c(G/N\times G/H, A)$ acts on $z\in C_c(G/N\times G\times G/H, A)\subseteq W$
by
\[
(b\cdot z)(rN,u,vH) = b(rN,vH)z(rN,u,vH).
\]
The group $G\times G/N$ is represented by
$(U\otimes\lambda^{G/N}\otimes 1)\times(1\otimes\rho\otimes\lambda^{G/H})$
and acts on $W$ by
\[
((s,tN)\cdot z)(rN,u,vH)
= \alpha_s(z(s^{-1}rtN,s^{-1}u,t^{-1}vH))
  \Delta_{G}(s)^{1/2}\Delta_{G/N}(tN)^{1/2}.
\]
Thus, for $\xi\in L^2(G/H,\H)$,
\begin{align*}
(\nu(b\cdot z)\xi)(rN,vH)
&= \int_G \pi\bigl((b\cdot z)(rN,y,vH)\bigr) U_y
\xi(y^{-1}rvH) \Delta_G(y)^{-1/2}\, dy\\
&= \int_G \pi\bigl(b(rN,vH)z(rN,y,vH)\bigr) U_y
\xi(y^{-1}rvH) \Delta_G(y)^{-1/2}\, dy\\
&= \pi(b(rN,vH)) \big((\nu(z)\xi)(rN,vH)\big)\\
&= \bigl( (\pi\otimes M^{G/N}\otimes
M^{G/H})(b)\bigr)(\nu(z)\xi)(rN,vH)\\
&=\nu_\llt(b)(\nu(z)\xi)(rN, vH),
\end{align*}
and
\begin{align*}
\big(\nu&((s,tN)\cdot z)\xi\big)(rN,vH)
= \int_G \pi\bigl((s,tN)\cdot z)(rN,y,vH)\bigr) U_y
\xi(y^{-1}rvH) \Delta_G(y)^{-1/2}\, dy\\
&= \int_G \pi\bigl(\alpha_s(z(s^{-1}rtN,s^{-1}y,t^{-1}vH))\bigr)
U_y \xi(y^{-1}rvH)
\Delta_{G}(s)^{1/2}\Delta_{G/N}(tN)^{1/2}\Delta_G(y)^{-1/2}\, dy\\
&= \int_G U_s\pi\bigl(z(s^{-1}rtN,s^{-1}y,t^{-1}vH)\bigr)
U_{s^{-1}y} \xi(y^{-1}rvH)
\Delta_G(s)^{1/2}\Delta_{G/N}(tN)^{1/2}\Delta_G(y)^{-1/2}\, dy\\
&\overset{(\dagger)}{=}
U_s\int_G \pi\bigl(z(s^{-1}rtN,y,t^{-1}vH)\bigr)
U_y \xi(y^{-1}s^{-1}rvH)
\Delta_{G/N}(tN)^{1/2}\Delta_G(y)^{-1/2}\, dy\\
&= U_s (\nu(z)\xi)(s^{-1}rtN,t^{-1}vH)
\Delta_{G/N}(tN)^{1/2}\\
&= \bigl((U_s\otimes\lambda_s\otimes 1)
(1\otimes\rho_{tN}\otimes\lambda^{G/H}_{tN})\bigr)
(\nu(z)\xi)(rN,vH)\\
&=\nu_\llt((s,tN))(\nu(z)\xi)(rN, vH),
\end{align*}
where the change of variables at $(\dagger)$ was $y\mapsto sy$.
Thus~(2) holds.

For (3), fix $\zeta>0$. 
Also fix nonzero $\phi\in C_c(G/N)$, $\eta\in C_c(G/H)$ and $h\in\H$.  It suffices to approximate (in $L^2(G/N\times G/H,\H)$) the function
\[
(rN, vH)\mapsto \phi(rN)\eta(vH)h.
\]
Using an approximate identity and the non-degeneracy of $\pi$, 
choose nonzero $a\in A$ such that $\|\pi(a)h-h\|<\zeta/(2\|\phi\otimes\eta\|_2)$.  Then choose a relatively compact open neighbourhood $O$ of $e$ in $G$ such that
$\|U_yh-h\|<\zeta/(2\|a\|\|\phi\otimes\eta\|_2)$ for $y\in O$.  Then
\[\|\pi(a)U_yh-h\|<\frac{\zeta}{\|\phi\otimes\eta\|_2}
\]
for all $y\in O$.  Next, choose $f\in C_c(G)$ with $\supp f\subseteq O$
such that $\int_G f(y)\Delta_G(y)^{-1/2}\, dy=1$. Also choose compact
subsets $L$ and $K$ of $G$ such that $L/N=\supp \phi$ and
$K/H=\supp\eta$, and $\xi\in C_c(G/H,\H)$ such that $\xi(s)=h$ for $sH\in O^{-1}LK/H$.  
Set
\[
z(rN,y,vH)=\phi(rN)f(y)\eta(vH)a.
\]
Then $z\in C_c(G/N\times G\times G/H, A)$ and
\begin{align*}
\|\nu(z)&\xi-\phi\otimes\eta\otimes h\|_2^2\\
&=\int_{G/N}\int_{G/H} \|\nu(z)\xi(rN, vH)-\phi(rN\eta(vH)h\|^2\, d(rN)\, d(vH)\\
&=\int_{G/N}\int_{G/H} \Big\|\int_G\phi(rN)\eta(vH)f(y)\pi(a)U_y\xi(y^{-1}rvH)\Delta_G(y)^{-1/2}\, dy
\\
&\hskip 4cm
-\Big(\int_G f(y)\Delta_G(y)^{-1/2}\, dy\Big)\phi(rN)\eta(vH)h\Big\|^2\, d(rN)\, d(vH)\\
&\leq \int_{G/N}\int_{G/H}\|\phi(rN)\eta(vH)\|^2\left(\int_G\|f(y)\Delta_G(y)^{-1/2}\|\|\pi(a)U_yh-h\|\, dy\right)^2\, d(rN)\, d(vH)
\end{align*}
by our choice of $\xi$.  Since
$\|\pi(a)U_yh-h\|<\zeta/(\|\phi\otimes\eta\|_2)$ for all $y\in\supp
f$ we have \[\|\nu(z)\xi-\phi\otimes\eta\otimes h\|_2<\zeta,\]
and hence $\nu$ is non-degenerate. Thus $(\nu_\llt,\nu,\nu_\rrt)$ is
a representation of $W$; the faithfulness follows because $\nu_\rrt$
is faithful.%
\footnote{Note that the proof that $(\nu_\llt,\nu,\nu_\rrt)$
is a faithful representation did not require $H$ to be normal in $G$,
but we do need the normality in what follows.}

We will obtain our representation of $K(\epsilon)\rtimes_{\epsilon_K|}
(G/H)$ by first constructing a representation $(\mu_\llt,\mu,\mu_\rrt)$
of  $K(\epsilon)=(A\rtimes_\alpha G)\otimes L^2(G/N)$ on the pair
$(\H\otimes L^2(G/N),\H)$ of Hilbert spaces, and then applying
Lemma~\ref{lem-coaction-x} to the coaction $\epsilon_{K}|$ of $G/H$
on $K(\epsilon)$.

We represent $A\rtimes_\alpha G\rtimes_\epsilon(G/N
)\rtimes_{\hat\epsilon}(G/N)$ on $\H\otimes L^2(G/N)$ and $A\rtimes_\alpha G$ on $\H$ by
\begin{align*}
\mu_\llt&:=(((\pi\rtimes U)\otimes\lambda)\circ\epsilon)\rtimes (1\otimes M^{G/N})\rtimes (1\otimes\rho)\\
&=(\pi\otimes 1)\rtimes (U\otimes\lambda^{G/N})\rtimes (1\otimes M^{G/N})\rtimes (1\otimes\rho)
\end{align*} and $\mu_\rrt:=\pi\rtimes U$,
respectively, and let
$
\mu:K(\epsilon)\to B(\H,\H\otimes L^2(G/N))
$
be the linear map such that
\[
\mu(b\otimes f)h=\pi\rtimes U(b)h\otimes f\quad\quad(b\in A\rtimes_\alpha G, f\in L^2(G/H), h\in\H).
\]
Note that $\mu$ is non-degenerate because $\pi\rtimes U$ is. For
$b\otimes f,c\otimes g\in K(\epsilon)$
and $h,k\in\H$,
\begin{align*}
\bigl( \mu(b\otimes f) h \bigm| \mu(c\otimes g)k\bigr)
&= \bigl( \pi\rtimes U(b)h \otimes f\bigm|\pi\rtimes U(c)k \otimes g \bigr)\\
&= \bigl(\pi\rtimes U(c^*b) h \bigm| k \bigr)(g\mid f)\\
&= \bigl(\pi\rtimes U(c^*b(g\mid f)) h \bigm| k \bigr)\\
&= \bigl( \pi\rtimes U(\langle c\otimes g\,,\, b\otimes f\rangle_{A\rtimes_\alpha G})h \bigm|k
\bigr) \\
&= \bigl( \mu_\rrt(\langle c\otimes g\,,\, b\otimes f\rangle_{A\rtimes_\alpha G}) h\bigm|k \bigr),
\end{align*}
so $\mu(b\otimes f)^*\mu(b\otimes f) =\mu_\rrt(\langle b\otimes f\,,\, b\otimes f\rangle_{A\rtimes_\alpha G})$.

On the pieces $A$, $C^*(G)$, $C_0(G/N)$ and  $C^*(G/N)$, $\mu_\llt$ is given  by
$\pi\otimes 1$, $U\otimes\lambda^{G/N}$, $1\otimes
M^{G/N}$ and  $1\otimes\rho$, respectively.
Let $(i_A,i_G)$ be the universal covariant representation of
$(A,G,\alpha)$. The left action of $A\rtimes_\alpha G\rtimes_\epsilon
(G/N) \rtimes_{\hat\epsilon}(G/N)$  on $K(\epsilon)$ is via the isomorphism
\[
((\id\otimes\lambda)\circ\epsilon)
\rtimes (1\otimes M^{G/N})
\rtimes (1\otimes\rho)
=(i_A\otimes 1)\rtimes (i_G\otimes\lambda^{G/N})\rtimes (1\otimes M^{G/N}     )
\rtimes (1\otimes\rho)
\]
of $A\rtimes_\alpha G \rtimes_\epsilon(G/N) \rtimes_{\hat\epsilon}
(G/N)$ onto $(A\rtimes_\alpha G)\otimes \K(L^2(G/N))$.
For $a\in A$, $g\in C_0(G/N)$ and $s\in G$ we have
\begin{align*}
\mu_\llt(a)\mu(b\otimes f)h&=(\pi(a)\otimes 1)(\pi\rtimes U(b)h\otimes f)\\
&=\pi\rtimes U(i_A(a)b)h\otimes f\\
&=\mu((i_A(a)\otimes1)(b\otimes f))h;\\
\mu_\llt(s)\mu(b\otimes f)h
&= (U_s\otimes\lambda_s^{G/N})(\pi\rtimes U(b)h\otimes f)\\
&=\pi\rtimes U(i_G(s)b)h\otimes\lambda^{G/N}_s(f)\\
&=\mu((i_G(s)\otimes\lambda^{G/N}_s)(b\otimes f))h;\\
\mu_\llt(g)\mu(b\otimes f)h&=(1\otimes M^{G/N}(h))(\pi\rtimes U(b)h\otimes f)\\
&=\pi\rtimes U(b)h\otimes M^{G/N}(h)f\\
&=\mu((1\otimes M^{G/N}(h))(b\otimes f)h;\\
\mu_\llt(sN)\mu(b\otimes f)h
&=1\otimes\rho_{sN}(\pi\rtimes U (b)h\otimes f)\\
&=\pi\rtimes U(b)h\otimes\rho_{sN}(f)\\
&=\mu(1\otimes\rho_{sN}(b\otimes f))h.
\end{align*}It follows that $(\mu_\llt, \mu,\mu_\rrt)$ is a representation of 
$K(\epsilon)$ on $(\H\otimes L^2(G/N),\H)$.

The coaction $\epsilon_K$ of $G/N$ on $K(\epsilon)$ is defined in \cite[Proposition~4.2]{KQ-MI} by
\[
\epsilon_K(b\otimes f)
= V(\epsilon(b)\otimes f)^{\Sigma_{23}}
\]
(here $\Sigma_{23}(b\otimes f\otimes g)=b\otimes g\otimes f$ and  $V
\in M((A\rtimes_\alpha G)\otimes \K(L^2(G/N))\otimes C^*(G/N))$
is given by
\[
V=1\otimes(M^{G/N}\otimes\id)(w_{G/N}^*),
\]
where $w_{G/N}\in UM(C_0(G/N)\otimes C^*(G/N))$ is the usual
multiplicative unitary).  
By Lemma~\ref{lem-coaction-x},
$\big(\mu_\llt\rtimes (G/H),\mu\rtimes (G/H),\mu_\rrt\rtimes (G/H)\big)$
is a faithful representation of $K\rtimes_{\epsilon|}(G/H)$  on the subspace
\begin{equation}
\label{Kset}
\clsp\{
((\mu\otimes\lambda)\circ\epsilon_K|(b\otimes f))(1\otimes M_g)
\mid b\in A\rtimes_\alpha G, f\in L^2(G/N), g\in C_0(G/H) \}
\end{equation}
of $B\big(\H\otimes L^2(G/H),\H\otimes L^2(G/N)\otimes
L^2(G/H)\big)$.  Since the canonical isomorphism of $(A\otimes
C_0(G/H))\rtimes_{\alpha\otimes\lt}G$ onto $A\rtimes_\alpha
G\rtimes_\epsilon(G/H)$ carries the representation $\nu_\rrt=(\pi\otimes
M^{G/H})\rtimes (U\otimes\lambda^{G/H})$ into 
$\mu_\rrt\rtimes (G/H)
=(((\pi\rtimes U)\otimes\lambda^{G/H})\circ\epsilon)
\rtimes (1\otimes M^{G/H})$,
the ranges of $\nu_\rrt$ and $\mu_\rrt\rtimes (G/H)$ clearly coincide.
(It is also not hard to check that the canonical isomorphism of the
left-hand coefficient algebras carries
$\nu_\llt$ into $\mu_\llt\rtimes(G/H)$.)

To finish the proof of the theorem, we need to show that the ranges  of
$\nu$ and $\mu\rtimes (G/H)$ coincide.  By the Rieffel correspondence it suffices to show that a dense subset of the range of $\nu$ is contained in the range of $\mu\rtimes (G/H)$.
To do this, we  need a more useful expression for
terms of the form
$(\mu\otimes\lambda)\circ\epsilon_K|(b\otimes f)$.
We have
\[(\mu\otimes\lambda)\circ\epsilon_K|
= (\mu\otimes\lambda\circ q)\circ\epsilon_K.\]
Recall that
$\epsilon= \hat\alpha|
= (\id\otimes q)\circ((i_A\otimes 1)\rtimes(i_G\otimes\iota))
= (i_A\otimes 1)\rtimes(i_G\otimes q)$,
where $\iota\colon G \to C^*(G)$ is the canonical map,
and $q$ maps $C^*(G)$ into $C^*(G/N)$.
Thus, for $b\in C_c(G,A)\subseteq A\rtimes_\alpha G$, we have
\[
\epsilon(b)
= \int_G (i_A\otimes 1)(b(s))(i_G\otimes q)(s) \, ds
= \int_G i_A(b(s))i_G(s) \otimes q(s) \, ds.
\]
Let $f\in L^2(G/N)$ and $h\in \H$. We write $c_f(h) = h\otimes f$ for $h\in\H$; note that
\[
\mu(b\otimes f) = c_f \circ (\pi\rtimes U(b)).
\]
We have
\[
(\epsilon(b)\otimes f)^{\Sigma_{23}}
= \int_G i_A(b(s))i_G(s) \otimes f\otimes q(s) \, ds,
\]
and therefore
\[
(\mu\otimes \lambda\circ q)(\epsilon(b)\otimes f)^{\Sigma_{23}}
= \int_G (c_f\circ \pi(b(s))U_s) \otimes \lambda_{sH}\, ds.
\]
In other words, for $\xi\in L^2(G/H,\H)$,
\[
(\mu\otimes \lambda\circ q)(\epsilon(b)\otimes f)^{\Sigma_{23}}
\xi(rN,tH)
= \int_G f(rN)\pi(b(s))U_s\xi(s^{-1}tH)\, ds.
\]

Next, for $T\in \K(L^2(G/N))$, $z,w\in C^*(G/N)$,
$f\in L^2(G/N)$, and $b\in A\rtimes_\alpha G$, compute:
\begin{align*}
(\mu\otimes\lambda\circ q)
((1\otimes T\otimes z)(b\otimes f\otimes w))
&= (\mu\otimes\lambda\circ q)
(b\otimes Tf \otimes zw)\\
&= c_{Tf}\circ(\pi\rtimes U(b))\otimes \lambda\circ q(zw)\\
&= \big((\id\otimes\id\otimes\lambda\circ q)(1\otimes T\otimes z)\big)
\big((\mu\otimes\lambda\circ q)(b\otimes f\otimes w)\big).
\end{align*}
From this we deduce that, for multipliers of $(A\rtimes_\alpha G )\otimes\K\otimes C^*(G/N)$ of the form $1\otimes m$, we have
\[
(\mu\otimes\lambda\circ q)(1\otimes m)
= 1\otimes (\id\otimes\lambda\circ q)(m),
\]
and in particular
\[
(\mu\otimes\lambda\circ q)(V)
= 1\otimes (M^{G/N}\otimes\lambda\circ q)(w_{G/N}^*);
\]
that is, for $\eta\in L^2(G/N\times G/H,\H)$,
\[
(\mu\otimes\lambda\circ q)(V)\eta(rN,vH)
= \eta(rN,rvH).
\]

Combining the above, we get
\begin{align*}
(\mu\otimes\lambda)\circ\epsilon_K|(b\otimes f)\xi(rN,vH)
&=(\mu\otimes \lambda\circ q)(\epsilon(b)\otimes f)^{\Sigma_{23}}
\xi(rN,rvH)\\
&= \int_G f(rN)\pi(b(y))U_y \xi(y^{-1}rvH))\, dy
\end{align*}
for $b\in C_c(G,A)\subseteq A\rtimes_\alpha G$,
$f\in L^2(G/N)$, and
$\xi\in L^2(G/H,\H)$.
The image of $K(\epsilon)\rtimes_{\epsilon_K|}(G/H)$ is thus
densely spanned by the operators defined by
\begin{equation*}
\label{mu}
(\mu\otimes\lambda)\circ\epsilon_K|(b\otimes f)((1\otimes M_g)\xi)(rN,vH)
= \int_G f(rN)\pi(b(y)) g(y^{-1}rvH) U_y \xi(y^{-1}rvH))\, dy.
\end{equation*}
Let $z\in C_c(G/N\times G\times G/H, A)$ be the function
$
z(rN,v,vH)=f(rN)b(y)g(y^{-1}rvH)\Delta_G(y)^{1/2};
$
then \begin{align*}((\mu\rtimes G/H)(b\otimes f)(1\otimes M_g))\xi(rN,vH)
&=(\mu\otimes\lambda)\circ\epsilon_K|(b\otimes f)(1\otimes M_g)\xi(rN,vH)\\
&=\int_G\pi(z(rN,y,vH))U_y\xi(y^{-1}rvH)\Delta_G(y)^{-1/2}\, dy\\
&=(\nu(z)\xi)(rN, vH).
\end{align*}
It follows that
the ranges of $\nu$ and $\mu\rtimes (G/H)$
in $B(\H\otimes L^2(G/H),\H\otimes L^2(G/N)\otimes L^2(G/H))$
coincide, and this completes the proof of Proposition~\ref{new-R}.
\end{proof}

%======================================================================

\section{Proof of Theorem~\ref{thm-KXZ}}
\label{next}

Recall that in Theorem~\ref{thm-KXZ} we assume that
$\alpha$ is a continuous action of a locally compact group $G$ 
by automorphisms of a $C^*$-algebra $A$, $N$ and $H$ are closed
normal subgroups of $G$ with $N\subseteq H$,  
and we have let 
$\epsilon$ denote the maximal coaction $\hat\alpha|_{G/N}$ of $G/N$
on $A\rtimes_\alpha G$.
We also retain the symmetric-imprimitivity bimodules
$W(P)$, $W(Q)$, and $W(R)$ defined in Section~\ref{previous},
and all the associated notation.

The basic idea is to invoke the symmetric imprimitivity
calculus of Theorem~\ref{thm-PhashQ}
and then show that $P\# Q$ is equivariantly isomorphic to
$R$, so that
\begin{align*}
X_{H/N}^{G/N}(\hat\epsilon)\otimes_{*} Z_{G/H}^{G/N}(\alpha)
&\cong W(P)\otimes_{\Phi}W(Q)
\cong W(P\# Q)
\cong W(R)
\cong K(\epsilon)\rtimes_{\epsilon_K|}(G/H).
\end{align*}
However, there are many isomorphisms of the coefficient algebras
involved here (see diagram~\eqref{biggie}), 
several of them non-canonical, and 
we must make sure they
are all compatible with this argument.

%----------------------------------------------------------------------
\subsection{Applying Theorem~\ref{thm-PhashQ} to $P$ and $Q$}

The map $(rN,u,vN)\mapsto
(v^{-1}rN,e)G= (v^{-1}N,r^{-1})G$ of $P$ into $Q/G$ induces a homeomorphism $\phi\colon K\under P\to Q/G$
such that
\[
\phi(K(rN,u,vN)) = (v^{-1}N,r^{-1})G.
\]
Moreover, $\phi$ is $L$-equivariant: for $(hN,y)\in L=N/H\times G$,
\begin{align*}
\phi(K(rN,u,vN)\cdot(hN,y))
&= \phi(K(ryN,uy,vhN))\\
&= (h^{-1}v^{-1}N,y^{-1}r^{-1})G\\
&= (h^{-1}N,y^{-1})\cdot (v^{-1}N,r^{-1})G\\
&= (hN,y)^{-1}\cdot \phi(K(rN,u,vN)).
\end{align*}
Thus the fibred product of $P$ and $Q$ over $\phi$ is
\begin{align*}
P\times_\phi Q
&= \{ (p,q)\in P\times Q \mid \phi(Kp)=qG\}\\
&= \{ (rN,u,vN,wN,z)\in G/N\times G\times G/N\times G/N\times G
\mid (v^{-1}N,r^{-1})G=(wN,z)G \}\\
&= \{ (rN,u,vN,wN,z)\in G/N\times G\times G/N\times G/N\times G
\mid wN=v^{-1}rzN \},
\end{align*}
and the right action of $L$ on $P\times_\phi Q$ is given by
\[
(rN,u,vN,wN,z)\cdot(hN,y) = (ryN,uy,vhN,h^{-1}wN,y^{-1}z).
\]
%Then $P\# Q := (P\times_\phi Q)/L$ has natural commuting free and
%proper actions of $K$ on the left and $G$ on the right, and we can
%apply the symmetric imprimitivity theorem to the data $({}_K(P\#Q)_G,A,
%\sigma,\tau)$.

Now define $\tilde\sigma\colon
P\to\Aut A$ and $\tilde\tau\colon Q\to\Aut A$ by
\[
\tilde\sigma_{(rN,u,vN)} = \alpha_u
\quad\text{and}\quad
\tilde\tau_{(wN,z)} = \id.
\]
We have
\[
\tilde\sigma_{(tN,s)\cdot(rN,u,vN)\cdot(hN,y)}
= \tilde\sigma_{(tryN,suy,tvhN)}
= \alpha_{suy}
= \alpha_s\alpha_u\alpha_y
= \sigma_{(tN,s)}\tilde\sigma_{(rN,u,vN)}
  \zeta_{(hN,y)}
\]
and
\[
\tilde\tau_{(hN,s)\cdot(wN,z)\cdot y}
= \tilde\tau_{(hwyN,szy)}
= \id
= \id\id\id
= \eta_{(hN,s)}\tilde\tau_{(wN,z)}\tau_{y}.
\]
It is clear that $\zeta$, $\sigma$, and $\tilde\sigma$ commute with
$\eta$, $\tilde\tau$, and $\tau$, since the latter are trivial.  Thus
all the hypotheses of Theorem~\ref{thm-PhashQ} are satisfied. 
Therefore, 
there exist isomorphisms
\begin{align*}
&\Phi\colon \Ind_K^P\sigma\rtimes_{\eta\otimes\rt} L
\to \Ind_G^Q\tau\rtimes_{\zeta\otimes\lt} L,\\
&\Phi_\sigma:\Ind_K^{P\#Q}\sigma\rtimes_{\tau\otimes\rt} G\to \Ind_L^Q\zeta\rtimes_{\tau\otimes\rt} G\text{\ and}\\ &\Phi_\tau:\Ind_G^{P\#Q}\rtimes_{\sigma\otimes\lt}K\to\Ind_L^P\eta\rtimes_{\sigma\otimes\lt}K
\end{align*}
such that the upper square of the following diagram commutes:
\begin{equation}\label{Phi}
\begin{split}
\xymatrix{
\Ind_G^{P\#Q}\tau\rtimes_{\sigma\otimes\lt}K
\ar[d]_{\Phi_\tau}
^\cong
\ar[rr]^{W(P\#Q)}
&& \Ind_K^{P\#Q}\sigma\rtimes_{\tau\otimes\rt}G
\ar[d]_{\Phi_\sigma}
^\cong
\\
\Ind_L^P\eta\rtimes_{\sigma\otimes\lt}K
\ar[rr]^{W(P)\otimes_\Phi W(Q)}
\ar[d]_{W(P)}
&& \Ind_L^Q\zeta\rtimes_{\tau\otimes\rt}G
\\
\Ind_K^P\sigma\rtimes_{\eta\otimes\rt}L
\ar[rr]^{\Phi}_\cong
&&
\Ind_G^Q\tau\rtimes_{\zeta\otimes\lt}L.
\ar[u]_{W(Q)}
}
\end{split}
\end{equation}
The lower square commutes by definition of $W(P)\otimes_\Phi W(Q)$.

Since we will need it later, we recall from Lemma~\ref{lem-one-and-two}
that the isomorphism $\Phi$ is induced by the $L$-equivariant isomorphism
$T\colon \Ind_K^P\sigma \to \Ind_G^Q\tau$ defined by
\begin{equation}
\label{eq-T}
T(f)(rN,s) = f(s^{-1}N,e,r^{-1}N)
\end{equation}
(because $\phi(K(s^{-1}N,e,r^{-1}N))=(rs^{-1}N,e)G=(rN,s)G$).
Further, the isomorphism $\Phi_\sigma$
is induced by the $G$-equivariant
isomorphism $\phi_\sigma\colon \Ind_K^{P\# Q}\sigma
\to \Ind_L^Q\zeta$ given by
\begin{equation}\label{eq-phisigma}
\phi_\sigma(f)(wN,z) = f((z^{-1}N,e,w^{-1}N,wN,z)L);
\end{equation}
$\Phi_\tau$ is induced by the $K$-equivariant
isomorphism $\phi_\tau\colon \Ind_G^{P\# Q}\tau
\to \Ind_L^P\eta$ given by
\begin{equation}\label{eq-phitau}
\phi_\tau(f)(rN,u,vN) = f((rN,u,vN,v^{-1}rN,e)L).
\end{equation}

%----------------------------------------------------------------------
\subsection{$P\#Q$ and $R$ are isomorphic}

The map $\psi\colon P\# Q \to R$ given by
\[
\psi((rN,u,vN,wN,z)L)  =  (rzN,uz,vH)
\]
is a (well-defined) homeomorphism with inverse
given by
$\psi^{-1}(rN,u,vH) = (rN,u,vN,v^{-1}rN,e)L$.
Since $\psi$ is equivariant for the left action
of $K$ and the right action of $G$, $\psi$ induces
induces a $K$-equivariant isomorphism
$\psi_\tau\colon \Ind_G^R\tau \to \Ind_G^{P\# Q}\tau$
such that
\begin{equation}\label{eq-psitau}
\psi_\tau(f)((rN,u,vN,wN,z)L) = f(rzN,uz,vH)
\end{equation}
and a $G$-equivariant isomorphism
$\psi_\sigma\colon \Ind_K^R\sigma \to \Ind_K^{P\# Q}\sigma$
with the same rule: 
\begin{equation}\label{eq-psisigma}
\psi_\sigma(f)((rN,u,vN,wN,z)L) = f(rzN,uz,vH).
\end{equation}

The map of $C_c(R,A)$ into $C_c(P\# Q,A)$
induced by $\psi$ extends
to an imprimitivity bimodule
isomorphism $\Psi\colon W(R)\to W(P\# Q)$
whose coefficient maps are
$\Psi_\tau:=\psi_\tau\rtimes K$
and $\Psi_\sigma:=\psi_\sigma\rtimes G$.
In other words, the following diagram commutes:
\begin{equation}
\label{R}
\xymatrix{
\Ind_G^{R}\tau\rtimes_{\sigma\otimes\lt}K
\ar[rr]^{W(R)}
\ar[d]_{\Psi_\tau=\psi_\tau\rtimes K}^\cong
&& \Ind_K^{R}\sigma\rtimes_{\tau\otimes\rt}G
\ar[d]^{\Psi_\sigma=\psi_\sigma\rtimes G}_\cong
\\
\Ind_G^{P\#Q}\tau\rtimes_{\sigma\otimes\lt}K
\ar[rr]^{W(P\#Q)}
&& \Ind_K^{P\#Q}\sigma\rtimes_{\tau\otimes\rt}G.
}
\end{equation}

%----------------------------------------------------------------------
\subsection{Assembly}

Now we assemble the commuting diagrams involving
the three bimodules from Theorem~\ref{thm-KXZ}
into diagram~\eqref{biggie} below.
(For simplicity we only indicate
the bimodules and isomorphisms, and the respective
diagram numbers.)
Note that every arrow is invertible, and 
the outer rectangle (whose vertical sides collapse) 
is precisely diagram~\eqref{eq-KXZ}.
Thus, to complete the proof of Theorem~\ref{thm-KXZ}, 
it remains to show that the squares labelled~\eqref{Lambda},
\eqref{Xi}, and~\eqref{T} commute, as well as the upper and lower
left-hand corners.

\begin{equation}\label{biggie}
\xymatrix{%::::::::::::::::::::::::::::::::::::::::::::::
{\cdot}
\ar[rrr]^{X_{H/N}^{G/N}(\hat\epsilon)}
\ar[d]_\cong^c
&&&
{\cdot}
\ar[rr]_\cong^d
\ar[d]_\cong^d
&&
{\cdot}
\ar[rrr]^{Z_{G/H}^{G/N}(\alpha)}
\ar[ddd]_{(\Theta\rtimes L)\circ i}
&&&
{\cdot}
\ar[ddd]^{\Omega\rtimes G}
\\%------------------------------------------------------
{\cdot}
\ar[dd]_{(\Gamma\rtimes K)\circ\iota}
\ar[rrr]_{X_{H/N}^{G/N}(\beta)}
&&&
{\cdot}
\ar[dd]^{(\Upsilon\rtimes L)\circ i}
&
\eqref{T}
&&
\eqref{eq-W(Q)}
&&
\\%------------------------------------------------------
&
\eqref{X}
&&&&&&&
\\%------------------------------------------------------
{\cdot}
\ar[rrr]_{W(P)}
&&&
{\cdot}
\ar[rr]_{\Phi}
&&
{\cdot}
\ar[rrr]_{W(Q)}
&&&
{\cdot}
\\%------------------------------------------------------
&&&&
\eqref{Phi}
&&&&
\\%------------------------------------------------------
&&
{\cdot}
\ar[uull]_{\Phi_\tau}
\ar[rrrr]^{W(P\#Q)}
&&&&
{\cdot}
\ar[uurr]^{\Phi_\sigma}
&&
\\%------------------------------------------------------
&
\eqref{Lambda}
&&&
\eqref{R}
&&&
\eqref{Xi}
&
\\%------------------------------------------------------
&&
{\cdot}
\ar[uu]^{\Psi_\tau}
\ar[rrrr]_{W(R)}
&&&&
{\cdot}
\ar[uu]_{\Psi_\sigma}
&&
\\%------------------------------------------------------
\\%------------------------------------------------------
{\cdot}
\ar[uuuuuu]^{(\Gamma\rtimes K)\circ\iota}
\ar[uurr]_{(\Lambda\rtimes K)\circ\iota}
&&&&
\eqref{KR}
&&&&
{\cdot}
\ar[uull]^{\Xi\rtimes G}
\ar[uuuuuu]_{\Omega\rtimes G}
\\%------------------------------------------------------
\\%------------------------------------------------------
{\cdot}
\ar[uu]_\cong^c
&&
{\cdot}
\ar[ll]_\cong^b
\ar[uull]_\cong^a
\ar[rrrr]_{K(\epsilon)\rtimes(G/H)}
&&&&
{\cdot}
\ar[uurr]_\cong
&&
}%:::::::::::::::::::::::::::::::::::::::::::::::::::::::
\end{equation}

%----------------------------------------------------------------------
\subsection{Non-canonical isomorphisms in~\eqref{biggie}}

%Next, we establish the relationships between all
%the non-canonical isomorphisms used above.

Using Equations~\eqref{eq-phitau}, \eqref{eq-psitau},
\eqref{eq-Lambda}, and~\eqref{eq-Gamma}, for any 
$f\in A\otimes C_0(G/N)\otimes C_0(G/H)$
and any $(rN,u,vN)\in R$ we have
\begin{align*}
\phi_\tau(\psi_\tau(\Lambda(f)))(rN,u,vN)
&= \psi_\tau(\Lambda(f))((rN,u,vN,v^{-1}rN,e)L)\\
& = \Lambda(f)(rN,u,vH)\\
& = f(ru^{-1}N,vH)\\
&= \Gamma(f)(rN,u,vN).
\end{align*}
So 
$\phi_\tau\circ\psi_\tau\circ\Lambda 
= \Gamma \colon A\otimes C_0(G/N)\otimes C_0(G/H)\to \Ind_L^R\eta$.
Since all four maps are $K$-equivariant,
it follows that the following diagram of isomorphisms commutes:
\begin{equation}
\label{Lambda}
\xymatrix{
\big((A\otimes C_0(G/N)\rtimes_{\alpha\otimes\lt} G)
\otimes C_0(G/H)\big)\rtimes_{\beta\otimes\lt}(G/N)
\ar[rr]^-{(\Lambda\rtimes K)\circ\iota}
%(A\otimes C_0(G/N)\otimes C_0(G/H))
%\rtimes_\epsilon(G/N\times G)
\ar[d]_{(\Gamma\rtimes K)\circ \iota}
&&\Ind_G^R\tau\rtimes_{\sigma\otimes\lt}K
\ar[d]_{\Psi_\tau=\psi_\tau\rtimes K}
\\
\Ind_L^P\eta\rtimes_{\sigma\otimes\lt}K
&&\Ind_G^{P\#Q}\tau\rtimes_{\sigma\otimes\lt}K.
\ar[ll]_-{\Phi_\tau = \phi_\tau\rtimes K}
}
\end{equation}

Using Equations~\eqref{eq-phisigma}, \eqref{eq-psisigma},
\eqref{eq-Xi}, and~\eqref{eq-Omega}, 
for any $f\in A\otimes C_0(G/H)$ and any $(wN,z)\in Q$ we have
\begin{align*}
\phi_\sigma(\psi_\sigma(\Xi(f)))(wN,z)
&= \psi_\sigma(\Xi(f))((N,z,w^{-1}N,wN,e)L)\\
&= \Xi(f)(N,z,w^{-1}H)\\
&= \alpha_z(f(w^{-1}H))\\
&= \Omega(f)(wN,z).
\end{align*}
Thus $\phi_\sigma\circ\psi_\sigma\circ\Xi = \Omega\colon
A\otimes C_0(G/H) \to \Ind_L^Q\zeta$.
All four maps are $G$-equivariant, so
the following diagram commutes:
\begin{equation}
\label{Xi}
\xymatrix{
\Ind_K^R\sigma\rtimes_{\tau\otimes\rt}G
\ar[d]_-{\Psi_\sigma=\psi_\sigma\rtimes G}
&&
(A\otimes C_0(G/H))\rtimes_{\alpha\otimes\lt}G
\ar[d]_-{\Omega\rtimes G}
\ar[ll]_-{\Xi\rtimes G}
\\
\Ind_K^{P\#Q}\sigma\rtimes_{\tau\otimes\rt}G
\ar[rr]^-{\Phi_\sigma= \phi_\sigma\rtimes G}
&&
\Ind_L^Q\zeta\rtimes_{\tau\otimes\rt}G.
}
\end{equation}

For $f\in A\otimes C_0(G/N)$ and any $(rN, s)\in G/N\times G$, 
using Equations~\eqref{eq-T}, \eqref{eq-Upsilon} and~\eqref{eq-Theta}, 
we have
\[
T(\Upsilon(f))(rN,s) = \Upsilon(f)(s^{-1}N,e,r^{-1}N)
= f(sr^{-1}N) = \Theta(f)(rN,s),
\]
so $T\circ\Upsilon = \Theta$.
All three maps are $L$-equivariant, so
the following diagram commutes:
\begin{equation}
\label{T}
\xymatrix{
\Ind_K^P\sigma\rtimes_{\eta\otimes\rt}L
\ar[rr]^{\Phi=T\rtimes L}
&&\Ind_G^Q\tau\rtimes_{\zeta\otimes\lt}L\\
&
\big((A\otimes C_0(G/N))\rtimes_{\alpha\otimes\lt}G\big)\rtimes_{\beta|}(H/N).
\ar[ul]^{(\Upsilon\rtimes L)\circ i}
\ar[ur]_{(\Theta\rtimes L)\circ i}
}
\end{equation}

%----------------------------------------------------------------------
\subsection{Canonical isomorphisms in~\eqref{biggie}}

For the upper left-hand square of diagram~\eqref{biggie},
temporarily set $C=A\rtimes_\alpha G \rtimes_\epsilon(G/N)$
and $D=(A\otimes C_0(G/N))\rtimes_{\alpha\otimes\lt}G$.
Then it is straightforward to verify that
the $\hat\epsilon-\beta$ equivariant canonical map
of $C$ onto $D$ induces an imprimitivity bimodule isomorphism
of $X_{H/N}^{G/N}(\hat\epsilon)$
onto $X_{H/N}^{G/N}(\beta)$
such that the diagram 
\begin{equation}\label{upper-left}
\xymatrix{%::::::::::::::::::::::::::::::::::::::::::::::
(C\otimes C_0(G/H))\rtimes_{\hat\epsilon\otimes\lt}(G/N)
%\ar[rr]^(.6){X_{H/N}^{G/N}(\hat\epsilon)}
\ar[rr]^-{X_{H/N}^{G/N}(\hat\epsilon)}
\ar[d]_\cong^c
&&
C\rtimes_{\hat\epsilon|}(H/N)
\ar[d]_\cong^d
\\%------------------------------------------------------
(D\otimes C_0(G/H))\rtimes_{\beta\otimes\lt}(G/N)
%\ar[rr]^(.6){X_{H/N}^{G/N}(\beta)}
\ar[rr]^-{X_{H/N}^{G/N}(\beta)}
&&
D\rtimes_{\beta|}(H/N)
}%:::::::::::::::::::::::::::::::::::::::::::::::::::::::
\end{equation}
commutes.
The lower left-hand triangle of~\eqref{biggie}, 
which is enlarged below, commutes because
all the isomorphisms are canonical.

\begin{equation*}
\xymatrix{%::::::::::::::::::::::::::::::::::::::::::::::
{(D\otimes C_0(G/H))\rtimes_{\beta\otimes\lt}(G/N)}
&
\\%------------------------------------------------------
\\%------------------------------------------------------
{(C\otimes C_0(G/H))\rtimes_{\hat\epsilon\otimes\lt}(G/N)}
\ar[uu]_\cong^c
&
{C\rtimes_{\hat\epsilon}(G/N)\rtimes_{\hat{\hat\epsilon}|}(G/H)}
\ar[l]_\cong^b
\ar[uul]_\cong^a
}%:::::::::::::::::::::::::::::::::::::::::::::::::::::::
\end{equation*}

This completes the proof of Theorem~\ref{thm-KXZ}.

%======================================================================

\section{Induction in Stages}
\label{main-proof}

We can deduce induction-in-stages for
the $Z$'s from results already in the literature:

\begin{prop}\label{prop-ZZZ}
Let $\alpha :G\to\Aut A$ be a continuous action of a locally compact group $G$ by automorphisms of a $C^*$-algebra $A$. Also let $H$ and $N$ be  closed
subgroups of $G$ with $N$ normal in $G$ and  $N\subseteq H$.
Then the following diagram of right-Hilbert bimodules commutes:
\[
\xymatrix{
(A\otimes C_0(G))\rtimes_{\alpha\otimes\lt}G
\ar[dr]_{Z_{G/N}^G(\alpha)}
\ar[rr]^{Z_{G/H}^G(\alpha)}
&
&
(A\otimes C_0(G/H))\rtimes_{\alpha\otimes\lt}G
\\
&
(A\otimes C_0(G/N))\rtimes_{\alpha\otimes\lt}G.
\ar[ur]_{Z_{G/H}^{G/N}(\alpha)}
&
\\
}
\]
\end{prop}

\begin{proof}
%As usual, we denote by $X_H^G(\alpha)$ Green's
%\[
%((A\otimes C_0(G/H))\rtimes_{\alpha\otimes\lt} G)
%-(A\rtimes_{\alpha|}H)
%\] 
%imprimitivity bimodule.
Consider the following diagram, where, as usual, 
the $X$'s denote Green imprimitivity bimodules:
%\begin{equation}
%\label{eq-dual}
%\begin{split}
%\xymatrix{
%(A\otimes C_0(G/N))\rtimes_{\alpha\otimes\lt}G
%\ar[ddd]_{Z_{G/H}^{G/N}(\alpha)}
%\ar[rr]^{X_N^G(\alpha)}
%&&
%A\rtimes_{\alpha|}N
%\ar[ddd]_{\Res}
%\\
%&
%((A\otimes C_0(G))\rtimes_{\alpha\otimes\lt}G)
%\rtimes_{\beta|}N
%\ar[d]_{\Res}
%\ar[ur]^{X_e^G(\alpha)\rtimes N}
%\ar[ul]_{Z_{G/N}^G(\alpha)}
%&
%\\
%&
%((A\otimes C_0(G))\rtimes_{\alpha\otimes\lt}G)
%\rtimes_{\beta|}H
%\ar[dl]^{Z_{G/H}^G(\alpha)}
%\ar[dr]_{X_e^G(\alpha)\rtimes H}
%&
%\\
%(A\otimes C_0(G/H))\rtimes_{\alpha\otimes\lt}G
%\ar[rr]_{X_H^G(\alpha)}
%&&
%A\rtimes_{\alpha|}H.
%}
%\end{split}
%\end{equation}
%\begin{equation}
%\label{eq-dual2}
%\begin{split}
%\xymatrix{
%A\rtimes_{\alpha|}N
%\ar[ddd]_{\Res}
%&&
%(A\otimes C_0(G))\rtimes_{\alpha\otimes\lt}G\rtimes_{\beta|}N
%\ar[ll]_{X_e^G(\alpha)\rtimes N}
%\ar[dl]_{Z_{G/N}^G(\alpha)}
%\ar[ddd]^{\Res}
%\\
%&
%(A\otimes C_0(G/N))\rtimes_{\alpha\otimes\lt}G
%\ar[ul]_{X_N^G(\alpha)}
%\ar[d]^{Z_{G/H}^{G/N}(\alpha)}
%&
%\\
%&
%(A\otimes C_0(G/H))\rtimes_{\alpha\otimes\lt}G
%\ar[dl]^{X_H^G(\alpha)}
%&
%\\
%A\rtimes_{\alpha|}H
%&&
%(A\otimes C_0(G))\rtimes_{\alpha\otimes\lt}G\rtimes_{\beta|}H.
%\ar[ll]^{X_e^G(\alpha)\rtimes H}
%\ar[ul]^{Z_{G/H}^G(\alpha)}
%}
%\end{split}
%\end{equation}
%--
%--
\begin{equation}\label{ZXX}
\xymatrix@!0@R=40pt@C=80pt{
&
&
A\rtimes_{\alpha|}N
\ar'[d][ddd]_{\Res}
&
\\
(A\otimes C_0(G))\rtimes_{\alpha\otimes\lt}G\rtimes_{\beta|}N
\ar[urr]^{X_e^G(\alpha)\rtimes N}
\ar[rrr]_{Z_{G/N}^G(\alpha)}
\ar[ddd]_{\Res}
& 
& 
&
(A\otimes C_0(G/N))\rtimes_{\alpha\otimes\lt}G
\ar[ul]_(.5){X_N^G(\alpha)}
\ar[ddd]^{Z_{G/H}^{G/N}(\alpha)}
\\
&&&\\
&
&
A\rtimes_{\alpha|}H
& 
\\
(A\otimes C_0(G))\rtimes_{\alpha\otimes\lt}G\rtimes_{\beta|}H
\ar[urr]^{X_e^G(\alpha)\rtimes H}
\ar[rrr]_{Z_{G/H}^G(\alpha)}
& 
& 
&
(A\otimes C_0(G/H))\rtimes_{\alpha\otimes\lt}G
\ar[ul]_(.5){X_H^G(\alpha)}
}
\end{equation}
%--
Commutativity of the right rear face
is exactly \cite[Proposition~3.5]{aHKRW-EP}.  
The upper and lower (triangular) faces
commute by \cite[Theorem~3.1]{EKR-CP}.%
\footnote{The statement of  \cite[Theorem~3.1]{EKR-CP} 
should end with ``-- $A\rtimes_{\alpha|}H$ bimodules''.}
The left rear face commutes by naturality of restriction
(\cite[Lemma~5.7]{KQR-DR}).
Since all except the vertical arrows are
imprimitivity bimodules, 
it follows that the front face commutes.

The commutative front face of diagram~\eqref{ZXX} should be viewed as a
strong version of induction in stages; 
the proposition follows from this because
\[
\Res\otimes_{( (A\otimes
C_0(G))\rtimes_{\alpha\otimes\lt}G)\rtimes_{\beta|}H} Z_{G/H}^G(\alpha)
\cong Z_{G/H}^G(\alpha)
\]
as a right-Hilbert $(A\otimes C_0(G))\rtimes_{\alpha\otimes\lt}G -
(A\otimes C_0(G/H))\rtimes_{\alpha\otimes\lt}G$ bimodule.
%The proposition follows quickly from this since the actions of $(A\otimes
%C_0(G))\rtimes_{\alpha\otimes\lt}G \rtimes_{\beta|}H $ and $(A\otimes
%C_0(G))\rtimes_{\alpha\otimes\lt}G \rtimes_{\beta|}N$ on
%$Z_{G/N}^G(\alpha)$
%and $\Res$, respectively, extend to the multiplier algebras.
\end{proof}

We next deduce induction-in-stages for the Mansfield bimodule in the
case of a dual coaction.
The hypotheses are the same as in Proposition~\ref{prop-ZZZ}.

\begin{prop}\label{Y-alpha}
Let $\alpha\colon G\to\Aut A$ 
be a continuous action of a locally compact group $G$ 
by automorphisms of a $C^*$-algebra $A$. 
Also let $H$ and $N$ be  closed
subgroups of $G$ with $N$ normal in $G$ and  $N\subseteq H$.
Then the following diagram of right-Hilbert bimodules commutes:
\[
\xymatrix{
A\rtimes_\alpha G \rtimes_{\hat\alpha}G
\ar[dr]_{Y_{G/N}^G(\hat\alpha)}
\ar[rr]^{Y_{G/H}^G(\hat\alpha)}
&
&
A\rtimes_\alpha G \rtimes_{\hat\alpha|}(G/H).
\\
&
A\rtimes G \rtimes_{\hat\alpha|}(G/N)
\ar[ur]_{Y_{G/H}^{G/N}(\hat\alpha|)}
&
\\
}
\]
\end{prop}

\begin{proof}
Consider the following diagram, where the 
isomorphisms are the canonical ones:%
%--

\begin{equation}\label{YtoZ2}
\xymatrix@!0@R=40pt@C=80pt{
(A\otimes C_0(G))\rtimes_{\alpha\otimes\lt}G\rtimes_{\beta|}N
\ar[rrr]^{Z_{G/N}^G(\alpha)}
\ar[ddd]_{\Res}
\ar[rd]^(.6){\cong}
&
&
&
(A\otimes C_0(G/N))\rtimes_{\alpha\otimes\lt}G
\ar'[d][ddd]^{Z_{G/H}^{G/N}(\alpha)}
\ar[rd]^(.6){\cong}
&
\\
&
(A\rtimes_\alpha G)\rtimes_{\hat\alpha}G\rtimes_{\hat{\hat\alpha}|}N
\ar[rrr]_{Y_{G/N}^G(\hat\alpha)}
\ar[ddd]^{\Res}
& 
& 
&
(A\rtimes_\alpha G)\rtimes_{\hat\alpha|}(G/N)
\ar[ddd]^{Y_{G/H}^{G/N}(\hat\alpha|)}
\\
&&&&\\
(A\otimes C_0(G))\rtimes_{\alpha\otimes\lt}G\rtimes_{\beta|}H
\ar'[r][rrr]^(.4){Z_{G/H}^G(\alpha)}
\ar[rd]_(.4){\cong}
&
&
&
(A\otimes C_0(G/H))\rtimes_{\alpha\otimes\lt}G
\ar[rd]_(.4){\cong}
& 
\\
&
(A\rtimes_\alpha G)\rtimes_{\hat\alpha}G\rtimes_{\hat{\hat\alpha}|}H
\ar[rrr]_{Y_{G/H}^G(\hat\alpha)}
& 
& 
&
(A\rtimes_\alpha G)\rtimes_{\hat\alpha|}(G/H).
}
\end{equation}
%--
The rear face is the commutative front face
of diagram~\eqref{ZXX}; 
the upper and lower faces commute by \cite[Proposition~1.1]{EKR-CP};
the right-hand face is seen to commute by ignoring the
left $H/N$-actions in Theorem~\ref{thm-YZ}; 
and 
it is straightforward to verify directly that the left-hand face commutes 
(or one can use naturality of restriction \cite[Lemma~5.7]{KQR-DR}).
It follows that the front face commutes, 
and the proposition follows from this
as in the proof of Proposition~\ref{prop-ZZZ}.
%Hence we have the following commutative
%diagram of right-Hilbert bimodules:
%\begin{equation}\label{strong-Ya}
%\xymatrix{
%(A\rtimes_\alpha G)\rtimes_{\hat\alpha}G \rtimes_{\hat{\hat\alpha}|}N
%\ar[rr]^-{Y_{G/N}^G(\hat\alpha)}
%\ar[d]_{\Res}
%&&
%(A\rtimes_\alpha G)\rtimes_{\hat\alpha|}(G/N)
%\ar[d]^{Y_{G/H}^{G/N}(\hat\alpha|)}
%\\
%(A\rtimes_\alpha G)\rtimes_{\hat\alpha}G \rtimes_{\hat{\hat\alpha}|}H
%\ar[rr]_-{Y_{G/H}^G(\hat\alpha)}
%&&
%(A\rtimes_\alpha G)\rtimes_{\hat\alpha|}(G/H).
%}
%\end{equation}
\end{proof}

%----------------------------------------------------------------------

\begin{proof}[Proof of Theorem~\ref{main-thm}]
Recall that we assume 
$\delta$ is a maximal coaction of a locally compact group $G$
on a $C^*$-algebra $B$, and that $N$ and $H$ are closed normal subgroups
of $G$ with $N\subseteq H$.  
Now let $(A,\alpha) = (B\rtimes_\delta G,\hat\delta)$,
and consider the following diagram:
\begin{equation}\label{Y2}
\xymatrix@!0@R=40pt@C=80pt{
(A\rtimes_\alpha G)\rtimes_{\hat\alpha}G\rtimes_{\hat{\hat\alpha}|}N
\ar[rrr]^{Y_{G/N}^G(\hat\alpha)}
\ar[ddd]_{\Res}
\ar[rd]^(.6){K(\delta)\rtimes G\rtimes N}
&
&
&
(A\rtimes_\alpha G)\rtimes_{\hat\alpha|}(G/N)
\ar'[d][ddd]^{Y_{G/H}^{G/N}(\hat\alpha|)}
\ar[rd]^(.6){K(\delta)\rtimes(G/N)}
&
\\
&
B\rtimes_{\delta}G\rtimes_{\hat{\delta}|}N
\ar[rrr]_{Y_{G/N}^G(\delta)}
\ar[ddd]^{\Res}
& 
& 
&
B\rtimes_{\delta|}(G/N)
\ar[ddd]^{Y_{G/H}^{G/N}(\delta|)}
\\
&&&&\\
(A\rtimes_\alpha G)\rtimes_{\hat\alpha}G\rtimes_{\hat{\hat\alpha}|}H
\ar'[r][rrr]^(.4){Y_{G/H}^G(\hat\alpha)}
\ar[rd]_(.4){K(\delta)\rtimes G\rtimes H}
&
&
&
(A\rtimes_\alpha G)\rtimes_{\hat\alpha|}(G/H)
\ar[rd]_(.4){K(\delta)\rtimes(G/H)}
& 
\\
&
B\rtimes_{\delta}G\rtimes_{\hat{\delta}|}H
\ar[rrr]_{Y_{G/H}^G(\delta)}
& 
& 
&
B\rtimes_{\delta|}(G/H).
}
\end{equation}
The rear face is the commutative front face
of diagram~\eqref{YtoZ2}; 
the upper, lower, and right-hand faces all
commute by naturality of the Mansfield bimodule
(\cite[Theorem~6.6]{KQ-MI}); 
the left-hand face commutes by naturality of restriction
(\cite[Lemma~5.7]{KQR-DR}).
The arrows connecting the rear face to the front face are
all imprimitivity bimodules, hence invertible;
it follows that the front face commutes, 
and the theorem follows from this
as in the proof of Proposition~\ref{prop-ZZZ}.
\end{proof}

\begin{rem}\label{KQ-defn}
The overall 
structure of our proof of Theorem~\ref{main-thm} has been:
using naturality to pass to dual coactions (diagram~\eqref{Y2}); 
in the dual case replacing Mansfield bimodules
by symmetric-imprimitivity bimodules (diagram~\eqref{YtoZ2}); 
and proving
induction-in-stages for symmetric-imprimitivity bimodules directly
(diagram~\eqref{ZXX}).

This amounts to gluing these three diagrams together along their common
faces, and in fact  we might have saved
some work by addressing the glued-together diagram directly
rather than the three separate pieces.
For example, 
part of the top face of the glued-together diagram
would be 
\begin{equation}\label{KQ}
\xymatrix{
(A\otimes C_0(G)) \rtimes_{\alpha\otimes\lt}G\rtimes_{\beta|}N
\ar[rr]^{Z_{G/N}^{G}(\alpha)}
\ar[d]_{\cong}
&&
(A\otimes C_0(G/N))\rtimes_{\alpha\otimes\lt}G
\ar[d]^{\cong}
\\
(A\rtimes_\alpha G)\rtimes_{\hat\alpha}G\rtimes_{\hat{\hat\alpha}|}N
\ar[d]_{K(\delta)\rtimes G\rtimes N}
\ar[rr]^{Y_{G/N}^{G}(\hat\alpha)}
&&
A\rtimes_\alpha G \rtimes_{\hat\alpha|}(G/N)
\ar[d]^{K(\delta)\rtimes(G/N)}
\\
B \rtimes_\delta G\rtimes_{\hat\delta|}N
\ar[rr]^{Y_{G/N}^{G}(\delta)}
&&
B \rtimes_{\delta|}(G/N),
}
\end{equation}
and the outer square of~\eqref{KQ}
is already known to commute: it is precisely the definition
of the Mansfield bimodule $Y_{G/N}^G(\delta)$
(\cite[Theorem~5.3]{KQ-MI}).  
While the argument may have been have made shorter in this way,
we feel that it is much 
better understood in terms of the three separate pieces.
\end{rem}

For future reference, we state 
as a corollary of the proof of Theorem~\ref{main-thm}
the strong version of induction in stages
which appears in diagram~\eqref{Y2}.
This is the analogue for maximal coactions  of
Theorem~4.1 of \cite{KQR-DR}, which was proved for
a (not-necessarily-maximal) 
coaction~$\delta$ of $G$ on $B$ and normal subgroups
$N\subseteq H$ of $G$ such that ``Mansfield imprimitivity works 
for~$H$''.

\begin{cor}\label{cor-KQR-analog}
Let $\delta: B\to M(B\otimes C^*(G))$ be a maximal coaction of a locally compact group $G$ on a $C^*$-algebra $B$.  Also let $N$ and $H$ be closed normal subgroups of $G$ such that $N\subseteq H$.  Then the diagram
\begin{equation}\label{eq-KQR-analog}
\xymatrix{
 B\rtimes_\delta G\rtimes_{\hat\delta|}N
\ar[rr]^-{Y_{G/N}^G(\delta)}
\ar[d]_{\Res}
&&
\quad B\rtimes_{\delta|}(G/N)\quad
\ar[d]^{Y_{G/H}^{G/N}(\delta|)}
\\
B\rtimes_\delta G\rtimes_{\hat\delta|}H
\ar[rr]_-{Y_{G/H}^G(\delta)}
&&
B\rtimes_{\delta|}(G/H)
}
\end{equation}
of right-Hilbert bimodules commutes.
\end{cor}

%======================================================================

\section{Another Application of Theorem~\ref{thm-PhashQ}}
\label{another}

Consider symmetric imprimitivity data $({}_KX_G,A,\sigma, \tau)$.  Then $({}_{\{e\}}X_G,A,\tau)$ is valid data as well, and  $W({}_{\{e\}}X_G)$ is an
$\Ind_G^X \tau $--$(C_0(X,A)\rtimes_{\tau\otimes\rt} G)$-imprimitivity bimodule which carries an action $\big(\sigma\otimes\lt,\sigma\otimes\lt,(\sigma\otimes\lt)\rtimes\id\big)$ of $K$.
Taking the crossed product of $W({}_{\{e\}}X_G)$ by the action of $K$
(see \cite{combes, cmw}) we get an
\[\big(\Ind_G^X \tau \rtimes_{\sigma\otimes\lt}K\big)-\big((C_0(X,A)\rtimes_{\tau\otimes\rt}
G)\rtimes_{\sigma\otimes\lt\rtimes\id}K\big)\] imprimitivity bimodule
$W({}_{\{e\}}X_G)\rtimes_{\sigma\otimes\lt}K$ which is a completion of
$C_c(K, C_c(X,A))$.  Similarly, $W({}_KX_{\{e\}})$ carries an action
$\big((\tau\otimes\rt)\rtimes\id,\tau\otimes\rt,\tau\otimes\rt\big)$ of
$G$, and taking crossed products by $G$ gives an
\[
\big((C_0(X,A)\times_{\sigma\otimes\lt}K)\rtimes_{\tau\otimes\rt\rtimes\id}
G\big)-\big(\Ind_K^X\sigma\rtimes_{\tau\otimes\rt}G\big)
\]
imprimitivity bimodule $W({}_KX_{\{e\}})\rtimes_{\tau\otimes\rt}G$. Let \[\Psi:(C_0(X,A)\rtimes_{\sigma\otimes\lt}K)\rtimes_{\tau\otimes\rt\rtimes\id}
G\to (C_0(X,A)\rtimes_{\tau\otimes\rt}
G)\rtimes_{\sigma\otimes\lt\rtimes\id}K\] be the natural isomorphism. It was proved in \cite[Lemma~4.8]{hrw} that there is an imprimitivity bimodule isomorphism
\begin{equation}
\label{eq-old-iso}
(W({}_{\{e\}}X_G)\rtimes K)\otimes_\Psi (W({}_KX_{\{e\}})\rtimes G)\cong W({}_KX_G),
\end{equation}
and it is an obvious test question for Theorem~\ref{thm-PhashQ} whether
it can recover this isomorphism on the level of spaces.

The first step is to note that
$W({}_{\{e\}}X_G)\rtimes_{\sigma\otimes\lt}K$ is isomorphic to the  imprimitivity
bimodule
$W({}_KP_{K\times G},A,\id,\sigma\times\tau)$
where $P:=K\times X$ and
\[
k\cdot (t,x)=(kt,x)\quad\text{and}\quad(t,x)\cdot (k,m)=(tk,k^{-1}\cdot x\cdot
m).
\]
To see this,  note that the map $K(t,x)\mapsto x$ is a homeomorphism of
$K\slash P$ onto $X$ and  $(t,x)(K\times G)\mapsto t\cdot xG$ is a homeomorphism
of $P/(G\times K)$ onto $X/G$, and
define
\begin{gather*}
\Lambda:\Ind_G^X \tau \to\Ind_{K\times G}^P(\sigma\times \tau)\text{\ \ by\ \ }
\Lambda(f)(t,x)=\sigma_t^{-1}(f(t\cdot x));\\
\Theta:C_0(X,A)\to \Ind_K^P \id\text{\ \ by\ \ }\Theta(h)(t,x)=h(x).
\end{gather*}
It is easy to check that $\Lambda$  and $\Theta$ are well-defined and
invertible, with inverses given by 
\[
\Lambda^{-1}(g)(x)=g(e,x)
\quad\text{and}\quad 
\Theta^{-1}(l)(e,x)=l(e,x)
\] 
for $g\in \Ind_{K\times G}^P(\sigma\times\tau)$ and  $l\in\Ind_K^P\id$.  
To check that $\Lambda$ is equivariant for the actions of $K$, it helps
to to write $\lt^X$ and $\lt^P$ to distinguish between actions induced
from left actions on different spaces. Then,
\begin{align*}
\Lambda\big((\sigma\otimes\lt^X)_k(f)\big)(t,x)
&=\sigma_t^{-1}\big((\sigma\otimes\lt^X)_k(f)(t\cdot x)\big )\\
&=\sigma_t^{-1}\sigma_k(f(k^{-1}t\cdot x))\\
&=\Lambda(f)(k^{-1}t,x)\\
&=\Lambda(f)(k^{-1}\cdot (t,x))\\
&=\big(\id\otimes\lt^P\big)_k(\Lambda(f))(t,x).
\end{align*}
Similarly, $\Theta$ is
$\big((\tau\times\sigma)\otimes(\rt^X\times\lt^X)\big)-\big((\sigma\times\tau)\otimes\rt^P\big)$
equivariant.
Thus $\Lambda$ and $\Theta$ induce  isomorphisms
\begin{gather*}\Lambda\rtimes
K:\Ind_G^X \tau\rtimes_{\sigma\otimes\lt^X}K\to\Ind_{K\times G}^P(\sigma\times
\tau)\rtimes_{\id\otimes\lt^P}K\\
\Theta\rtimes (G\times
K):C_0(X,A)\rtimes_{(\tau\times\sigma)\otimes(\rt^X\times\lt^X)}(G\times K)\to
\Ind_K^P \id \rtimes_{(\sigma\times\tau)\otimes\rt^P}(K\times G).
\end{gather*}
For $z\in C_c(K,C_c(X,A))$ define \[\Upsilon(z)(t,x)=\sigma_t^{-1}
\big(z(t)(t\cdot x)\big)\Delta_K(t)^{1/2}.\]  It is not hard to check, using the
formulas given  at
\cite[Equations~(B.2)]{BE} for the symmetric imprimitivity theorem bimodules and
at
\cite[Equations~3.5--3.8]{hrw-proper} for  the Combes crossed
product, that
$(\Lambda\rtimes K,\Upsilon,\Theta\rtimes (G\times K))$ extends to an
imprimitivity bimodule isomorphism of $W({}_{\{e\}}X_G)\times_{\sigma\otimes\lt}K$
onto $W({}_KP_{K\times G})$.

Similarly, $W({}_KX_{\{e\}})\rtimes_{\tau\otimes\rt}G$ is isomorphic to the  imprimitivity
bimodule associated to the data $({}_{K\times G}Q_G,A,\sigma\times\tau,\id)$ where $Q:=G\times X$ and
\[
(s,x)\cdot m=(m^{-1}s,x)\quad\text{and}\quad(k,m)\cdot(s,x)=(sm^{-1},k\cdot
x\cdot m^{-1}).
\]
(In place of $(\Lambda,\Upsilon,\Theta)$  use $(\Gamma,\Omega,\xXi)$ where, for $s\in G$,
\begin{gather*}
\Gamma:\Ind_K^X\sigma\to \Ind_{K\times G}^Q(\sigma\times\tau)\text{\ \ is\ \ }\Gamma(f)(s,x)=\tau_s^{-1}(f(x\cdot s^{-1}));\\
\xXi: C_0(X,A)\to\Ind_G^Q\id\text{\ \ is\ \ }\xXi(h)(s,x)=h(x);\text{\ and}\\
\Omega:C_c(G, C_c(X,A))\to W(Q)\text{\ \ is\ \ }\Omega(z)(s,x)=z(s^{-1},x)\Delta_G(s)^{-1/2}.)
\end{gather*}

The hypotheses of Theorem~\ref{thm-PhashQ} are satisfied with $\phi:K\under P\to
Q/G$ given by
$\phi(K(t,x))=(e,x)G$ and  $\tilde\sigma_{(t,x)}=\sigma_t$ and
$\tilde\tau_{(s,x)}=\tau_s$.
Thus
\begin{align*}
P\times_\phi Q&=\{ (t,x,s,y):t\in K, s\in G, x,y\in X \text{\ and\ }
\phi(K(t,x))=(s,y)G\}\\
&=\{(t,x,s,x):t\in K, s\in G,  x\in X\}
\end{align*}
and $K\times G$ acts on $P\times_\phi Q$ by the diagonal action
\[
(t,x,s,x)\cdot (k,m)=(tk, k^{-1}\cdot x\cdot m, sm, k^{-1}\cdot x\cdot m).
\]
The map $\psi:P\times_\phi Q\to X$ given by $(t,x,s,x)\mapsto t\cdot x\cdot s^{-1}$
induces a homeomorphism $\bar\psi$ of $P\# Q=(P\times_\phi Q)/(K\times G)$ onto $X$.  Then $\bar\psi$ is equivariant for the actions of $K$ and $G$ because $\psi$ is:  for $k\in K$ and $m\in G$ we have
\begin{gather*}
k\cdot\psi(t,x,s,x)=k\cdot(t\cdot x\cdot s^{-1})=kt\cdot x\cdot
s^{-1}=\psi(kt,x,s,x)=\psi(k\cdot (t,x,s,x))\\
\psi(t,x,s,x)\cdot m=(t\cdot x\cdot s^{-1})\cdot m=t\cdot x\cdot
s^{-1}m=\psi(t,x,m^{-1}s)=\psi((t,x,s)\cdot m).
\end{gather*}
Thus $W({}_K(P\# Q)_G)$ and $W({}_KX_G)$ are isomorphic. The isomorphism \eqref{eq-old-iso} now follows from Theorem~\ref{thm-PhashQ}.

%======================================================================

\end{document}